\documentclass{amsart}
\usepackage{graphicx}
\usepackage{amssymb}
\usepackage{epstopdf}
\usepackage{nicefrac}
\usepackage{enumerate}
\usepackage{hyperref}
\usepackage{float}
\usepackage{color}
\usepackage{pst-all}
\usepackage{tabularx}

\newtheorem{thm}{THEOREM}[section]

\newtheorem{cor}[thm]{COROLLARY}
\newtheorem{defn}[thm]{DEFINITION}

\newtheorem{lemma}[thm]{LEMMA}

\newtheorem{prop}[thm]{PROPOSITION}

\newcommand{\ds}{\displaystyle}

\newcommand{\D}{{{\mathbb D}^q}}

\newcommand{\ab}{\alpha\beta}

\newcommand{\dF}{d_{\F}}

\newcommand{\diam}{\mbox{\rm diam\,}}

\newcommand{\vol}{\mbox{\rm vol\,}}

\newcommand{\e}{{\epsilon}}
\newcommand{\ez}{{\epsilon_0}}
\newcommand{\eone}{{\epsilon_1}}
\newcommand{\etwo}{{\epsilon_2}}
\newcommand{\ethree}{{\epsilon_3}}
\newcommand{\estar}{{\epsilon_*}}

 \newcommand{\fF}{\mathfrak{F}}

\newcommand{\hh}{{\bf h}}

\newcommand{\F}{{\mathcal F}}

\newcommand{\cPa}{{\cP}_{\alpha}}

\newcommand{\wcP}{\widetilde{\cP}}

\newcommand{\cTa}{{\cT}_{\alpha}}

\newcommand{\cTab}{{\cT}_{\alpha \beta}}
\newcommand{\cTba}{{\cT}_{\beta \alpha}}
\newcommand{\wcT}{\widetilde{\cT}}

\newcommand{\wcTab}{\widetilde{\cT}_{\alpha \beta}}
\newcommand{\wcTba}{\widetilde{\cT}_{\beta \alpha}}
\newcommand{\hab}{{\bf h}_{\alpha \beta}}

\newcommand{\whab}{\widetilde{\bf h}_{\alpha \beta}}

\newcommand{\HF}{{\cH}_{\F}}

\newcommand{\mR}{{\mathbb R}}

\newcommand{\cA}{{\mathcal A}}

\newcommand{\cH}{{\mathcal H}}
\newcommand{\cI}{{\mathcal I}}

\newcommand{\cN}{{\mathcal N}}

\newcommand{\cP}{{\mathcal P}}
\newcommand{\cQ}{{\mathcal Q}}

\newcommand{\cT}{{\mathcal T}}
\newcommand{\cU}{{\mathcal U}}

\newcommand{\cZ}{{\mathcal Z}}

\newcommand{\wtH}{{\widetilde H}}

\newcommand{\wtL}{{\widetilde L}}
\newcommand{\wtM}{{\widetilde M}}

\newcommand{\wtU}{{\widetilde U}}

\newcommand{\wtF}{{\widetilde{\mathcal F}}}

\newcommand{\wtx}{{\widetilde{x}}}

\newcommand{\wtpi}{{\widetilde{\pi}}}

\newcommand{\wtg}{{\widetilde{\gamma}}}
\newcommand{\wtvol}{{\widetilde{vol}}}

\def\subtrans{\mathbin{\cap{\mkern-9mu}\mid}\,\,}

\newcommand{\cat}{\mbox{\rm cat}}
\newcommand{\catt}{\cat_{\subtrans}}

\begin{document}

\title{Compact foliations with finite transverse LS category}

\begin{abstract}
We  prove that  if  $F$ is a   foliation of a compact manifold $M$ with all
leaves compact submanifolds,  and  the    transverse  saturated  category
of $F$ is finite,  then the leaf space $M/F$ is compact Hausdorff.
The proof is surprisingly delicate, and is based on some new observations
about  the geometry of compact foliations. Colman proved in her 1998 doctoral thesis
  that the transverse  saturated  category of a compact
Hausdorff foliation is always finite, so we obtain a new characterization of
the compact Hausdorff foliations among the compact foliations as those with
  finite transverse  saturated  category.
\end{abstract}

 \author{Steven Hurder}
 \author{Pawe\l\/ G. Walczak}
 \email{hurder@uic.edu, pawelwal@math.uni.lodz.pl}
\address{SH: Department of Mathematics, University of Illinois at Chicago,
322 SEO (m/c 249), 851 S. Morgan Street, Chicago, IL 60607-7045}
\address{PGW: Wydzia\l  ~Matematyki i Informatyki, Uniwersytet ~ \L\'odzki, Banacha 22, 90-238, \L\'od\'z, Poland}

\thanks{Version date: December 9, 2016}

\thanks{2010 {\it Mathematics Subject Classification}. Primary: 57R30, 53C12, 55M30; Secondary 57S15}

\thanks{Keywords: compact   foliation, transverse Lusternik-Schnirelmann category,  Epstein filtration}

\maketitle

\section{Introduction}\label{sec-intro}

A {\it compact foliation} is a foliation of a manifold $M$  with all leaves
compact submanifolds. For codimension one or two,
a compact foliation $\F$ of a compact manifold $M$  defines a fibration of
$M$ over  its leaf space $M/\F$ which is an orbifold
\cite{Reeb1952,Epstein1972,Epstein1976,Vogt1976,EMS1977}.
For codimension three and above,  the leaf space $M/\F$ of a compact
foliation need not be a Hausdorff space. This was first shown  by an example
of Sullivan \cite{Sullivan1976} of a flow on a 5-manifold whose orbits are
circles, and the lengths of the orbits are not bounded above.    Subsequent
examples of Epstein and Vogt
\cite{EV1978,Vogt1977b}   showed that for any codimension greater than two,
there are examples of  compact foliations of compact manifolds whose leaf
spaces are not  Hausdorff.
Vogt  gave a remarkable example   of a $1$-dimensional, compact
$C^0$-foliation of $\mR^3$ with no upper bound on the lengths of the circle
leaves in   \cite{Vogt1989}.   The results described below    apply to the
case of compact $C^1$-foliations of compact manifolds.

A compact foliation $\F$ with Hausdorff leaf space is said to be {\it
compact Hausdorff}. The holonomy of each leaf of a compact Hausdorff
foliation is a finite group, a property which characterizes them among the
compact foliations. If every leaf has trivial holonomy group, then a compact
Hausdorff foliation is a fibration. Otherwise,  a compact Hausdorff
foliation is  a   ``generalized Seifert fibration'', where the leaf space
$M/\F$ is a ``V-manifold''  \cite{Satake1956,Holmann1964,Millett1974}.

 A  compact foliation whose leaf space is non-Hausdorff  has  a
closed, non-empty saturated subset, the  {\it bad set} $X_1 \subset M$,
which is the union of the leaves whose holonomy group is infinite. The image
of $X_1$  in the leaf space $M/\F$ consists of the non-Hausdorff points for
the $T_1$ topology on $M/\F$.
The work by Edwards, Millett and Sullivan \cite{EMS1977}  established many
fundamental  properties of the geometry of the leaves of a compact
foliation near its bad set, yet  there is no   general structure theory for
compact foliations, comparable to what is understood for compact Hausdorff
foliations. The   results of \S\S4, 5 and 6 of this work provide new
insights and techniques for the study of these foliations. In particular, we
introduce the notion of a \emph{tame point} for the bad set $X_1$ in
Definition~\ref{def-pp}. This notion  is a new property of non-Hausdorff
compact foliations, and is  studied in Section~\ref{sec-tame}.

The transverse     Lusternik-Schnirelmann (LS) category of  foliations was
introduced in the 1998 thesis of   H. Colman    \cite{Colman1998,CM2000}.
The key idea is that of  a transversally categorical open set.
Let $(M,{\F})$ and $(M',{\F}')$ be foliated     manifolds. A
 homotopy $H \colon M' \times [0,1] \to M$ is said to be {\it
foliated}
if for all  $0 \leq t \leq 1$ the map $H_t$ sends each leaf
$L'$ of $\F'$ into another leaf $L$ of ${\F}$.
An open subset $U$ of $M$ is {\it transversely categorical} if  there
is a foliated  homotopy $H \colon U \times [0,1] \to M$ such that
$H_0 \colon U \to M$ is the  inclusion, and  $H_1\colon U \to M$ has image
in a
single leaf of $\F$. Here $U$ is regarded as a foliated  manifold with the
foliation induced by $\F$ on $U$.
In other words, an open subset
$U$ of $M$ is  transversely categorical if the inclusion
$(U,{\F}_U)\hookrightarrow (M,\F)$
factors through a single leaf, up to  foliated homotopy.

\begin{defn}
The  transverse (saturated) category $\catt  (M,\F)$ of a foliated  manifold
$(M,\F)$ is the least number of transversely categorical open saturated sets
required to cover $M$. If no such finite covering exists, then
$\catt  (M,\F)=\infty$.
\end{defn}

The   transverse    category $\catt(M,\F)$ of a  compact Hausdorff foliation
$\F$  of a compact manifold $M$ is  always  finite  \cite{CM2000}, as every
leaf   admits a saturated product neighborhood which is transversely
categorical.  For a non-Hausdorff compact foliation,   our main result is
that there is  no  transversely categorical covering of the   the bad set.
 \begin{thm}\label{thm-compact}
Let $\F$ be a compact  $C^1$-foliation of a  compact manifold $M$ with
non-empty bad set $X_1$.
Then there exists a dense set of  tame points  $X_1^t \subset X_1$.
Moreover, for each $x \in X_1^t$, there is no transversely categorical
saturated open set containing $x$.
\end{thm}

 \begin{cor}\label{cor-compact}
Let $\F$ be a compact  $C^1$-foliation of a   compact    manifold $M$. If
$M$ admits a covering by transversely categorical open saturated sets,  then
$\F$ is compact   Hausdorff.
Hence, the   transverse holonomy groups of all leaves of $\F$ are finite
groups, and $\F$ admits a transverse Riemannian metric which is holonomy
invariant.
\end{cor}
\proof
  Millett \cite{Millett1974} and Epstein \cite{Epstein1976}  showed that for
a  compact Hausdorff foliation $\F$ of a manifold $M$ each leaf has finite
holonomy, and thus  $M$ admits a Riemannian metric so that the foliation is
Riemannian.
\endproof

 Recall that a foliation is {\it geometrically taut} if the manifold $M$
admits a Riemannian metric so that each leaf is an immersed minimal manifold
\cite{Rummler1979,Sullivan1979,Haefliger1981}.  Rummler  proved in
\cite{Rummler1979} that a compact foliation is Hausdorff if and only if it
is taut, and thus we can conclude:
\begin{cor}\label{cor-taut}
A compact $C^1$-foliation of a   compact  manifold $M$ with  $\catt
(M,{\F}) < \infty$
is   geometrically taut.
\end{cor}

The idea of the proof of Theorem~\ref{thm-compact} is as follows.  The
formal definition of the bad set $X_1$ in \S 3 is that it consist of leaves
of $\F$ such that every open neighborhood of the leaf contains leaves of
arbitrarily large volume.
This   characterization of the bad set intuitively suggests that it should
be a rigid set.  That is,  any foliated homotopy of an  open neighborhood
of a point in the  bad set   should preserve  these  dynamical properties,
hence the open neighborhood  cannot be continuously retracted to a single
leaf.  The  proof of this   statement  is surprisingly delicate, and
requires a very precise understanding of the properties of leaves in an open
neighborhood of the bad set. A key result is   Proposition~\ref{prop-tame},
an extension of the Moving Leaf Lemma in  \cite{EMS1977},
 which establishes the existence of ``tame points''.

 The overview of the paper is as follows:
The first two sections consist of background material, which we recall to
establish notations, and also present a variety of technical  results
required in the later sections.
 In \S2 we give some basic results from foliation theory, and  in \S3 we
recall some basic results about compact foliations, especially the structure
theory for the good and the  bad sets.
In \S4 we   establish a key homological  property for compact leaves under
deformation by a homotopy. The techniques introduced in this section are
used again in later sections.
The most technical results of the paper are contained in   \S5, where  we
prove that tame points are dense in the bad set.
  Finally,  in \S6 we prove that an open saturated set containing a tame
point is not  categorical.
 Theorem~\ref{thm-compact} follows immediately from
Propositions~\ref{prop-tame} and \ref{prop-npp}.

 \section{Foliation preliminaries}\label{sec-prelims}

We assume that $M$ is a compact smooth   Riemannian manifold  without
boundary   of dimension $m = p +q$, that $\F$  is a compact  $C^1$-foliation
of codimension-$q$, and that the leaves of $\F$ are  smoothly immersed
compact submanifolds, so that $\F$ is more precisely   a
$C^{1,\infty}$-foliation.
  For   $x \in M$,   denote by $L_x$ the leaf of $\F$ containing $x$.

We recall below some well-known facts about foliations, and introduce some
conventions of notation. The books
\cite{CandelConlon2000, Godbillon1991, HectorHirsch1986} provide excellent
basic references; our notation is   closest to that used in
\cite{CandelConlon2000}.
Note that the analysis of the bad sets in later sections  requires careful
estimates  on  the foliation geometry;
 not just in each leaf, but also for nearby leaves of a given leaf.  This
requires a careful description of the  local metric geometry of a foliation,
as given in this section.

 \subsection{Tangential  and normal  geometry}\label{subsec-geodesic}

Let $T\F$ denote the tangent bundle to $\F$,   and let $\Pi \colon Q \to M$
denote its normal bundle, identified with the subbundle   $T\F^{\perp}
\subset TM$ of vectors orthogonal to $T\F$. The Riemannian metric on $TM$
induces  Riemannian metrics on both $T\F$ and $Q$  by fiberwise restriction.
For a vector $\vec{v} \in T_xM$, let $\|\vec{v}\|$ denote its length in the
Riemannian metric. Then for  $\vec{v} \in T_x\F$  the length in the induced
leafwise metric is also denoted by $\|\vec{v}\|$.

For $\e > 0$, let $T^{\e}M  \subset TM$ denote the disk subbundle  of
vectors with length less than $\e$, and let $T^{\e}{\F}  \subset T\F$ and
$Q^{\e} \subset Q$ be the corresponding    $\e$-disk subbundles of $T\F$ and
$Q$, respectively.

Let $d_M \colon M \times M \to [0,\infty)$ be the distance function
associated to the Riemannian metric on $M$.
Given  $r > 0$ and a set $K \subset  M$,  let
\begin{equation}
B_M(K,r)   =    \{ y \in M \mid d_M(K,y) < r\} \ .
\end{equation}

For a leaf $L \subset M$, let $d_L \colon L \times L \to [0,\infty)$ be the
distance function on  $L$ for the restricted Riemannian metric on  $L$. That
is,
  for $x,x' \in L$
the distance $d_L(x,x')$ is the infimum of the  lengths of piece-smooth
leafwise paths between $x$ and $x'$.
As $L$ is compact, the manifold  $L$ with  the   metric  $d_L$ is a complete
metric space, and the distance $d_L(x,x')$ is realized by a leafwise
geodesic path from $x$ to $x'$.
 We introduce the notation
$\dF$ for the collection of leafwise distance functions, where $\dF(x,y) =
d_L(x,y)$ if $x,y \in L$,  and otherwise $\dF(x,y) = \infty$.
Given  $r > 0$ and a set  $K \subset L$,  let
\begin{equation}\label{eq-leafwisenbhd}
B_{\F}(K,r)   =    \{ y \in M \mid \dF(K,y) < r\} \subset L \ .
\end{equation}

Let $\exp = \exp^{M} \colon TM \to   M$ denote the    exponential map
for $d_M$ which is well-defined as $M$ is compact. For
    $x \in M$, we let $\exp^M_x   \colon T_x M \to M$ denote the
exponential map at $x$.

For $x \in L$, we let $\exp^{\F}_x \colon T_xL \to L$ denote the
exponential map for the leafwise Riemannian metric.
Then  $\exp^{\F}_x$ maps the ball $B_{T_xL}(0,r)$ of radius $r$ in $T_x L$
onto the set  $B_{\F}(x,r)$.

We next chose $\ez > 0$ so that it satisfies a sequence of conditions, as
follows.
 For each $x \in M$, the    differential $D_{\vec{0}} \exp^M_x \colon T_xM
\cong T_{\vec{0}}(T_xM) \to T_xM$ is   the identity map. It follows that
there exists $\e_x > 0$ such that the restriction
 $\exp^M_x \colon T^{\e_x}_x M \to M$ is a diffeomorphism.  As $M$ is
compact,
there exists $\ez > 0$ such that for all $x \in M$,    the restriction
$\exp^M_x \colon T^{\ez}_xM  \to M$ is a diffeomorphism onto its image.
Thus,
$\ez$ is less than the injectivity radius of the Riemannian metric on $M$.
(See \cite{BC1964,doCarmo1992}  for details of the properties of the
injectivity radius of the geodesic map.)

We also require that $\ez >0$ be chosen so that  for all $x \in M$:
\begin{enumerate}
\item     The open ball
  $B_M(x,\ez)$ is  a totally normal neighborhood of $x$ for the metric $d_M$
This means that for any pair of points
$y,z \in B_M(x,\ez)$ there is a unique geodesic  contained  in $B_M(x,\ez)$
between $y$ and $z$.
In particular, $B_M(x,\ez)$ is geodesically convex  (See \cite[page
72]{doCarmo1992}.)
\item   The  leafwise exponential map    $\exp^{\F}_x \colon  T^{\ez}_x \F
\to L_x $ is a diffeomorphism onto it image.
\item   $B_{\F}(x,\ez) \subset L_x$  is a totally normal neighborhood of $x$
for the leafwise metric  $\dF$.
\end{enumerate}

Let $\exp^Q_x \colon Q_x \to M$ denote the restriction of  $\exp^M_x$ to the
normal bundle at $x$.
Then for all $x \in M$,  $\exp^Q_x \colon Q_x^{\ez} \to M$ is a
diffeomorphism onto its image.
We also require that $\ez > 0$ satisfy:
\begin{enumerate}\setcounter{enumi}{3}
\item For all $x \in M$, $\exp^Q_x \colon Q_x^{\ez} \to M$ is transverse to
$\F$, and  that the
image
 $\exp^Q_x (Q_x^{\ez})$ of the normal disk  has angle at least $\pi/4$ with
the leaves of the foliation $\F$.
\end{enumerate}

  We use the normal exponential map to define a normal product neighborhood
of a subset $K \subset L$ for a leaf $L$.
  Given $0< \e \leq \ez$, let  $Q(K,\e) \to K$ denote  the restriction of
the $\e$-disk bundle $Q^{\e} \to M$ to  $K$.
 The  normal neighborhood $\cN(K,\e)$ is the image of the  map,    $\exp^Q
\colon Q(K, \e) \to M$.
If $K = \{x\}$ is a point and $0 < \e < \ez$, then $\cN(x,\e)$ is a
uniformly transverse  normal   disk  to $\F$.

The restriction of the ambient metric   $d_M$ to a leaf $L$ need not
coincide (locally) with the leafwise geodesic metric $\dF$ -- unless the
leaves of $\F$ are totally geodesic submanifolds of $M$. In any case, the
Gauss Lemma implies that the two metrics are locally equivalent. We
require that $\ez > 0$ satisfy:
\begin{enumerate}\setcounter{enumi}{4}
\item For all $x \in M$, and
 for all   $y,y' \in B_{\F}(x,\ez)$,  then $\dF$ and $d_M$ are related
by
\begin{equation}  \label{eq-uniform}
 d_M(y,y')/2 ~ \leq ~ \dF(y,y') ~ \leq ~ 2 \, d_M(y,y') \ .
 \end{equation}
 \end{enumerate}

Let $dvol$ denote the leafwise volume $p$-form associated to the Riemannian
metric on $T\F$.
Given any   bounded, Borel     subset $A \subset L$ for the leafwise metric,
define its leafwise volume by $  vol(A) = \int_A \; dvol$.

 Let $L \subset M$ be a compact leaf, then there exists $0 < \e_L < \ez$
such that the normal geodesic map  $\exp^Q \colon Q(K, \e_L) \to M$
 is a diffeomorphism onto the open neighborhood $\cN(K,\e_L)$. We thus
obtain a normal projection map along the normal geodesic balls to points in
$L$,
 which we denote by $\Pi_L  \colon \cN(K,\e_L) \to L$. Note that the
restriction of $\Pi_L$ to $L$ is the identity map.

 Let $\F | \cN(K,\e_L)$ denote the restricted foliation whose leaves are the
connected components of the leaves of $\F$ intersected with $\cN(K,\e_L)$.
The tangent bundle to $\F | \cN(K,\e_L)$ is just the restriction of $T\F$ to
$\cN(K,\e_L)$, so for   $x' \in \cN(K,\e_L)$ and $x = \Pi_L(x')$,
 the differential of $\Pi_L$ induces a linear isomorphism $D_{\F}\Pi_L
\colon  T_{x'}\F \to T_xL$.
 Then the assumption on $\ez$ in \ref{subsec-geodesic}(4) above, implies
that $D_{\F}\Pi_L$ satisfies a Lipschitz estimate for some constant $C$,
which is the identity  when
restricted to the the leaf tangent bundle   $TL$.

 We use this observation in two ways.
For $L \subset M$ a compact leaf, assume that $0 < \e_L \leq \ez$ satisfies,
for $x \in L$ and  $x' \in \cN(K,\e_L)$ such that $x = \Pi_L(x')$:

  \begin{enumerate}\setcounter{enumi}{5}
\item  for the leafwise Riemannian volume $p$-form $dvol_{\F}$
\begin{equation}\label{eq-univol}
(dvol_{\F}|_{x'})/2 ~ \leq ~  (D_{\F}\Pi_L)^*(dvol_{\F}|_x)  ~  \leq ~ 2 \,
(dvol_{\F}|_{x'}) \ ;
 \end{equation}
\item  for the leafwise Riemannian norm $\| \cdot \|_{\F}$ and $\vec{v'} \in
T_{x'}\F$,
\begin{equation}\label{eq-unimetric}
(\|\vec{v'}\|_{\F})/2 ~ \leq ~  \|D_{\F}\Pi_L(\vec{v'})\|_{\F} ~  \leq ~ 2
\,  (\|\vec{v'}\|_{\F}) \ .
 \end{equation}
 \end{enumerate}

\subsection{Regular Foliation Atlas}\label{subsec-atlas}

We next recall some basic properties of foliation charts.
A {\it regular foliation atlas}  for $\F$ is a finite collection  $\{
(U_{\alpha},\phi_{\alpha})  \mid    \alpha \in \cA \}$ so that:
\begin{enumerate}
\item[(F1)] $\cU = \{ U_{\alpha}  \mid  \alpha \in \cA  \}$  is a covering
of $M$ by $C^{1,\infty}$--coordinate charts $\phi_{\alpha} : U_{\alpha}
\rightarrow  (-1,1)^m $  where each  $U_{\alpha}$ is a convex subset    with
respect to the metric $d_M$.
\item[(F2)]  Each  coordinate  chart $\phi_{\alpha} : U_{\alpha} \rightarrow
(-1,1)^m $
admits an extension to a $C^{1,\infty}$--coordinate chart  \\
$\widetilde{\phi}_{\alpha} : \wtU_{\alpha} \rightarrow (-2,2)^m$ where
$\wtU_{\alpha}$  is a convex subset    containing the $2\ez$-neighborhood of
$U_{\alpha}$, so $B_M(U_{\alpha}, \ez) \subset \wtU_{\alpha}$. In
particular, the closure $\overline{U_{\alpha}}  \subset \wtU_{\alpha}$.
\item[(F3)]  For each $z \in (-2,2)^q$, the preimage
$ \wcP_{\alpha}(z) = \widetilde{\phi}_{\alpha}^{-1} ( (-2,2)^{p} \times
\{z\}) \subset \wtU_{\alpha}$
 is  the connected component containing
$\widetilde{\phi}_{\alpha}^{-1}(\{0\} \times \{z\})$ of the intersection of
the leaf of $\F$ through $\phi_{\alpha}^{-1}(\{0\} \times \{z\})$ with the
set $\wtU_{\alpha}$.
\item[(F4)]   $\cP_{\alpha}(z)$ and  $\wcP_{\alpha}(z)$ are convex subsets
of diameter less than $1$
with respect to   $\dF$.
\end{enumerate}
The construction of regular coverings is described  in   chapter  1.2  of
\cite{CandelConlon2000}.

If the tangent bundle $T\F$ and normal bundle $Q = T\F^{\perp}$ to $\F$ are
oriented, then we   assume that the charts in the regular covering preserve
these orientations.

The inverse images
\[ \cP_{\alpha}(z) =   {\phi}_{\alpha}^{-1} ( (-1,1)^{p} \times \{z\})
\subset  {U}_{\alpha}\]
are smoothly embedded discs contained in the leaves of $\F$,   called the
{\em plaques} associated to the given foliation atlas.
The convexity hypotheses in (F4) implies  that
if $U_{\alpha}  \cap  U_{\beta} \not= \emptyset$, then each plaque
$\cP_{\alpha}(z)$ intersects at most one plaque of  $U_{\beta}$. The
analogous statement holds for pairs $\wtU_{\alpha}  \cap  \wtU_{\beta} \not=
\emptyset$.  More generally,  an  intersection of plaques
$\cP_{\alpha_1}(z_1) \cap \cdots \cap \cP_{\alpha_d}(z_d) $ is either empty,
or a convex set.

Recall that a   Lebesgue number for the covering $\cU$ is a constant $\e >
0$ so that    for each $x \in M$ there exists $U \in \cU$ with $B_M(x, \e)
\subset U$.  Every covering of a compact Riemannian manifold (in fact, of a
compact metric space) admits a Lebesgue number.   We also require that $\ez
> 0$ satisfy:
\begin{enumerate}\setcounter{enumi}{7}
\item $2 \ez$ is a Lebesgue number for the   covering  $\{  U_{\alpha}  \mid
\alpha \in \cA \}$ of $M$ by foliation charts.
 \end{enumerate}
  Then for any $x \in M$, the restriction of $\F$ to  $B_M(x,\ez)$  is a
product foliation, and by condition (F1) the   leaves  of $\F \mid
B_M(x,\ez)$ are convex discs for the metric $\dF$.

For each $\alpha \in \cA$, the extended chart $\widetilde{\phi}_{\alpha}$
defines a
$C^1$--embedding
$$ \widetilde{t}_{\alpha} = \widetilde{\phi}_{\alpha}^{-1}(\{0\} \times
\cdot) : (-2,2)^q
\rightarrow \wtU_{\alpha} \subset M
$$
whose image is denoted by $\wcT_{\alpha}$.  We can   assume that the  images
$\wcT_{\alpha}$ are pairwise disjoint. 
Let $t_{\alpha}$ denote the restriction of $\widetilde{t}_{\alpha}$ to
$(-1,1)^q \subset (-2,2)^q$, and
define  $\ds {\cT}_{\alpha} = t_{\alpha}(-1,1)^q$.
Then the collection of all plaques for the foliation atlas is  indexed by
the
{\em complete transversal}
$$ \cT = \bigcup_{\alpha \in {\cA}} \cT_{\alpha} \ . $$
For a point $x \in \cT$,   let
$\cP_{\alpha}(x) = \cP_{\alpha}(t_{\alpha}^{-1}(x))$ denote the plaque
containing $x$.

The Riemannian metric on $M$ induces a Riemannian metric and corresponding
distance function
${\bf d}_{\cT} $ on each extended  transversal $\wcT_{\alpha}$. For $\alpha
\not= \beta$ and $x \in \cT_{\alpha}$,  $y \in \cT_{\beta}$ we set ${\bf
d}_{\cT}(x,y) = \infty$.

Given  $x \in \wcT_{\alpha}$ and  $r > 0 $, let
$\ds {\bf B}_{\cT}(x,r)   =  \{ y \in \wcT_{\alpha}  \mid   {\bf
d}_{\cT}(x,y) < r \} $.

Given a subset $\cZ \subset U_{\alpha}$ let $\cZ_{\cP}$ denote the union of
all plaques in $U_{\alpha}$ having non-empty  intersection with  $\cZ$.
We set $\cZ_{\cT} = \cZ_{\cP} \cap \cT_{\alpha}$.  If $\cZ$ is an open
subset of $U_{\alpha}$, then $\cZ_{\cP}$ is open in  $U_{\alpha}$ and
$\cZ_{\cT}$ is an open subset of $\cT_{\alpha}$.

Given any point $w \in (-1,1)^p$,  we can define a transversal
 $\ds {\cT}_{\alpha}(w) = \phi_{\alpha}^{-1}(\{w\} \times (-1,1))$. There is
a canonical map
 $\ds \psi_{w} \colon  {\cT}_{\alpha}(w)  \to  {\cT}_{\alpha}(0) =
{\cT}_{\alpha}$ defined by, for $y \in (-1,1)^q$,
 \begin{equation} \label{eq-trans}
  \psi_w ( \phi_{\alpha}^{-1}(w \times \{y\}) =  \phi_{\alpha}^{-1}(0 \times
\{y\}) \ .
  \end{equation}

The Riemannian metric on $M$ induces also induces a Riemannian metric and
distance function on each transversal
$\ds {\cT}_{\alpha}(w)$. By mild abuse of notation we denote all such
transverse metrics by ${\bf d}_{\cT}$.
Then by the uniform extension property of the foliation charts,  there
exists a constant $C_T \geq 1$ so that
for all $\alpha \in \cA$,  $w \in (-1,1)^p$ and $x,y  \in
{\cT}_{\alpha}(w)$,
\begin{equation}\label{eq-CT1}
  {\bf d}_{\cT}(x,y)/C_T ~  \leq ~  {\bf d}_{\cT}(\psi_{w}(x),\psi_{w}(y)) ~
\leq  ~ C_T\, {\bf d}_{\cT}(x,y) \ .
\end{equation}
We use the maps \eqref{eq-trans}  to translate points in the coordinate
charts $U_{\alpha}$ to the ``center'' transversal $\cT_{\alpha}$.
The constant $C_T$ is a uniform estimate of the normal distortion introduced
by this translation.

We will  also consider the normal geodesic $\e$-disk $\cN(y,\e)$ at    $y =
\phi_{\alpha}^{-1}(w \times \vec{0})$, defined as the image of the map
$\ds \exp^Q_y \colon Q_y^{\e} \to \cN(y,\e)$,
which for $0 < \e \leq \ez$ is uniformly  transverse to $\F$.
path-length distance function on $\cN(y,\e)$ for the restricted Riemannian
metric on this transverse disk is equal to the distance function $d_M$ on
$M$ restricted to $\cN(y,\e)$.

Assume that  the image $\cN(y,\e) \subset U_{\alpha}$,  then we can project
it to the transversal $\cTa$ along the plaques in $U_{\alpha}$. Denote this
projection by
$\ds \Pi_{\alpha}^{\F} \colon \cN(y,\e) \to \cTa$.
We also assume that the constant    $C_T \geq 1$ is sufficiently large so
that for all $y \in M$, for all $0 < \e \leq \ez$, for all $\alpha$ with
$\cN(y,\e) \subset U_{\alpha}$ and for all
$z,z' \in \cN(y,\e)$ we have
\begin{equation}\label{eq-CT2}
  {\bf d}_{M}(z , z')/C_T ~  \leq ~  {\bf
d}_{\cT}(\Pi_{\alpha}^{\F}(z),\Pi_{\alpha}^{\F}(z')) ~ \leq  ~ C_T\, {\bf
d}_{M}(z,z') \ .
\end{equation}

 \subsection{Transverse holonomy}\label{subsec-holonomy}

The main result of this section is the definition of the    module of
uniform continuity function   for    elements of $\HF^n$, and its
application in Lemma~\ref{lem-ext}.

We  first recall the definition  of the holonomy pseudogroup of $\F$.
A pair of indices    $(\alpha, \beta)$  is said to be {\it admissible} if
$U_{\alpha} \cap U_{\beta} \not= \emptyset$.
Let $\cTab \subset \cTa$ denote  the open set of plaques of $U_{\alpha}$
which intersect some plaque of $U_{\beta}$. The   holonomy transformation
$\hab \colon \cTab \to \cTba$ is defined by
$y = \hab(x)$ if and only if $\cP_{\alpha}(x)  \cap \cP_{\beta}(y)  \not=
\emptyset$.
The finite collection
\begin{equation}\label{eq-holonomy}
\HF^1 = \{\hab \colon \cTab \to \cTba \mid  (\alpha, \beta)  ~ {\rm
admissible} \}  \ .
\end{equation}
generates the holonomy pseudogroup $\HF$ of local homeomorphisms of   $\cT$.

A {\it plaque chain of length $n$}, denoted by   ${\cP}$,  is  a collection
of   plaques
$$\{\cP_{\alpha_0}(z_0), \cP_{\alpha_1}(z_1), \ldots , \cP_{\alpha_n}(z_n)
\} $$
 satisfying $\cP_{\alpha_i}(z_i) \cap \cP_{\alpha_{i+1}}(z_{i+1}) \not=
\emptyset$ for $0 \leq i < n$.
Each pair of indices $(\alpha_i , \alpha_{i+1})$ is admissible, so
determines a    holonomy map
$\hh_{\alpha_i \alpha_{i+1}}$   such that
$\hh_{\alpha_i \alpha_{i+1}}(z_i) = z_{i+1}$.
Let $\hh_{\cP}$ denote the composition of these maps, so that
$$\hh_{\cP} = \hh_{\alpha_{n-1} \alpha_n} \circ \cdots \circ \hh_{\alpha_1
\alpha_2} \circ \hh_{\alpha_0 \alpha_1} \ . $$

Let  $\HF^n = \{\hh_{\cP} \mid \cP ~ {\rm has \ length \ at \ most \ n}\}
\subset   \HF$ denote the collection of maps obtained from the composition
of at most $n$ maps in  $\HF^1$.

 Each generator $\hab \colon \cTab \to \cTba$ is the restriction of the
transition map $\whab \colon \wcTab \to \wcTba$ defined by the intersection
of the extended charts
 $\widetilde{U}_{\alpha} \cap \widetilde{U}_{\beta}$. The domain $\cTab
\subset \wcTab$ is precompact with
${\bf B}_{\wcT}(\cTab,\ez) \subset \wcTab$, so $\hab$ is a uniformly
continuous   homeomorphism on its domain.
That is,    given any $0 < \e < \ez$, there   is a module of continuity
$\mu_{\ab}(\e) > 0$  such that for all   $x \in \cTab$
$${\bf B}_{\wcT}(x,\mu_{\ab}(\e)) \subset \wcTab \quad {\rm and} \quad
\whab({\bf B}_{\wcT}(x,\mu_{\ab}(\e))) \subset {\bf B}_{\wcT}(\hab(x),\e) \
. $$
For the admissible pairs $(\alpha, \alpha)$ we set   $\mu_{\alpha
\alpha}(\e) = \e$.
 Given $0 < \e \leq \ez$,  define
\begin{equation} \label{def-mu}
\mu(\e)  =  \min \{\mu_{\ab}(\e) \mid  (\alpha, \beta)  ~ {\rm admissible}
\}
\end{equation}
so that $0 < \mu(\e) \leq \e$.  Then for every admissible pair $(\alpha,
\beta)$ and each  $x \in \cTab$ the holonomy map $\hab$ admits an extension
to a   local homeomorphism $\whab$ defined by the holonomy of $\F$, which
satisfies
$\ds  \whab({\bf B}_{\wcT}(x,\mu(\e))) \subset {\bf B}_{\wcT}(\hab(x),\e)$.

For an integer $n > 0$ and $0 < \e \leq \ez$  recursively define
$\mu^{(1)}(\e) = \mu(\e)$
and $\mu^{(n)}(\e) = \mu(\mu^{(n-1)}(\e))$, so that $\mu^{(n)}$ denotes the
$n$-fold composition.
Then define
\begin{equation}
 \mu(n, \e) =  \min \{\e, \mu(\e), \mu( \mu(\e)),   \ldots, \mu^{(n)}(\e)\}
 \label{eq-unifcont}
 \end{equation}
Note that $0 < \mu(\e) \leq \e$ implies  $ \mu(n, \e) \leq \mu^{(n)}(\e)
\leq \e$.

 \begin{lemma} \label{lem-ext}
Given a   plaque chain $\cP$ of length $n$, and $0 < \e \leq \ez$ set
$\delta = \mu(n, \e)$.
 Then for any     $x$ in the domain of $\hh_{\cP}$,
     there  is an extension to a    local homeomorphism  $\widetilde
h_{\cP}$ defined by the holonomy of $\F$   whose domain includes the closure
of the  disk ${\bf B}_{\wcT}(x,\delta)$ about $x$ in $\wcT$, and
\begin{equation}\label{eq-unifn}
  \widetilde h_{\cP}({\bf B}_{\wcT}(x,\delta)) \subset {\bf
B}_{\wcT}(\hh_{\cP}(x),\e)   \ .
  \end{equation}
That is,   $ \mu(n, \e)$ is a   module of uniform continuity  for  all
elements of $\HF^n$.
 \end{lemma}
 \proof  For each $0 \leq i < n$,  $\mu(n, \e) \leq \mu(i, \e)$ hence  there
is an extension of
 $$\hh_i = \hh_{\alpha_{i-1} \alpha_i} \circ \cdots \circ \hh_{\alpha_0
\alpha_1}$$
to $\widetilde{\hh_i}$   whose domain includes the disk ${\bf
B}_{\wcT}(x,\delta)$ about $x$.
The image $\hh_i({\bf B}_{\wcT}(x,\delta)$ is contained in a ball of radius
at most $\mu(n-1, \e)$,  so that we can continue the extension process to
$\hh_{i+1}$.
\endproof

 \subsection{Plaque length and metric geometry}\label{subsec-metric}

We make two   observations about the metric leafwise geometry of foliations
\cite{Plante1975}.  In particular, the technical result
Proposition~\ref{prop-volume} below is a key fact for our proof of the main
result of this work.

Let  $\gamma \colon [0,1] \to L$ be a leafwise $C^1$-path. Its
  leafwise Riemannian length is  denoted  by  $|| \gamma ||_{\F}$.

The {\it plaque length}    of $\gamma$, denoted by   $|| \gamma ||_{\cP}$,
is the least integer $n$ such that the
image of $\gamma$ is covered by a chain of convex plaques
$$\{\cP_{\alpha_0}(z_0), \cP_{\alpha_1}(z_1), \ldots , \cP_{\alpha_n}(z_n)
\} $$
where $\gamma(0) \in \cP_{\alpha_0}(z_0)$, $\gamma(1) \in
\cP_{\alpha_1}(z_1)$, and  successive plaques $\cP_{\alpha_i}(z_i) \cap
\cP_{\alpha_{i+1}}(z_{i+1}) \not= \emptyset$.

\begin{prop}\label{prop-length}
For any leafwise $C^1$-path $\gamma$,
$\ds    ||\gamma ||_{\cP} ~  \leq ~  \lceil (||\gamma ||_{\F}/\ez) \rceil$.
Moreover, if $\gamma$ is leafwise geodesic, then
$\ds   ||\gamma ||_{\F} ~ \leq ~  ||\gamma ||_{\cP} + 1$.
\end{prop}
\proof   Let $N = \lceil (||\gamma ||_{\F}/\ez) \rceil$ be the least integer
greater than $||\gamma ||_{\F}/\ez$.  Then there exist points $0=  t_0 < t_1
< \cdots  < t_{N} = 1$
such that
the restriction of $\gamma$ to each segment $[t_i, t_{i+1}]$ has length at
most $\ez$.
The diameter of the set $\gamma([t_i, t_{i+1}])$ is at most $\ez$, hence
there is some
$U_{\alpha_i} \in \cU$ with $\gamma([t_i, t_{i+1}]) \subset U_{\alpha_i}$
hence $\gamma([t_i, t_{i+1}]) \subset {\cPa}(z_i)$ for some $z_i$.
Thus,  the image of $\gamma$ is covered by a chain of convex plaques of
length at most  $N$.

Conversely, suppose $\gamma$ is a leafwise geodesic  and  $\ds
\{\cP_{\alpha_0}(z_0), \cP_{\alpha_1}(z_1), \ldots , \cP_{\alpha_n}(z_n) \}
$ is a plaque chain covering the image $\gamma([0,1])$.
 Each plaque $\cP_{\alpha_i}(z_i)$ is a leafwise convex set of diameter at
most $1$ by the assumption (F4) in Section~\ref{subsec-atlas}, so
$\ds  ||\gamma ||_{\F} ~ \leq ~ (n+1) \leq ~  ||\gamma ||_{\cP} + 1$.
\endproof

The extension property (F2)  in Section~\ref{subsec-atlas} implies that  for
all $ \alpha \in \cA$ and  $z \in (-1,1)^q$,  the closure
$ \overline{{\cPa}(z)}$ is compact, hence has finite leafwise volume which
is uniformly continuous with respect to the parameter $z$.
 Hence,  there exist  constants  $0 < C_{min} \leq C_{max}$ such that
 \begin{equation}\label{eq-CP}
C_{min} ~ \leq ~ vol(\cP_{\alpha}(z))~  \leq ~ C_{max}, ~ \forall \alpha \in
\cA, ~ \forall z \in [-1,1]^{q} \ .
\end{equation}
We note a consequence of this uniformity which is critical to the proof of
the main theorem.

\begin{prop}\label{prop-volume}
Let $M$ be a compact manifold. Then there exists a monotone increasing
function $v \colon [0,\infty) \to [0,\infty)$ such that if $L$ is a compact
leaf,    then $vol(L) \leq v(\diam(L))$.
Conversely,   there exists a monotone increasing function $R \colon
[0,\infty) \to [0,\infty)$ such that if $L$ is a compact leaf,    then
$\diam(L) \leq R(vol(L))$.
\end{prop}
\proof

The holonomy pseudogroup of $\F$ has a finite set of generators, hence has
a \emph{uniform} upper bound on the   growth rate of words. This implies
that given $r > 0$,
there exists a positive integer $e(r)$ such that any subset of a leaf with
leaf diameter at most $r$ can be covered by no more that $e(r)$ plaques.
Thus, if $L$ is a leaf with diameter at most $r$, then $L$      has volume
at most $v(r) = C_{max} \cdot e(r)$.

Now suppose that  $L$ is a compact leaf with   diameter $r = \diam(L)$.
Then, for any pair of points $x,y \in L$, there exists a
 length   minimizing geodesic segment $\gamma \colon [0,1] \to L$ of length
$r = \dF(x,y)$, with $\gamma(0) = x$ and $\gamma(1) = y$.
Let $B_{\F}(\gamma,\ez)$ denote the leafwise $\ez$-tubular neighborhood of
the image of $\gamma$, defined by \eqref{eq-leafwisenbhd}.
Recall that   the restricted metric $\dF$ on leaves has uniformly bounded
geometry.  Then as $\ez$ is assumed in assumption (2) of
Section~\ref{subsec-geodesic} to be less than the injectivity radius for the
leafwise metric,  and $\gamma$ is  a length-minimizing geodesic, there is a
constant $V_0 > 0$ so that the leafwise volume $\vol(B_{\F}(\gamma,\ez))
\geq V_0 \cdot r$.
 Thus, $vol(L) > V_0 \cdot \diam(L)$, and then set $R(v) = v/V_0$.
 \endproof

 \subsection{Captured leaves} \label{subsec-captured}
  The main result is   Proposition~\ref{prop-cptcov}, which shows that given
a compact leaf $L$ of $\F$, and another  compact leaf $L'$   which is
sufficiently close to $L$ at some point, where how close depends  on
$vol(L')$,   then $L'$ is ``captured'' by the holonomy of $L$. We require
some preliminary definitions and observations before giving the proof of
this key fact.

Let $L \subset  M$ be a compact  leaf, and recall that in
Section~\ref{subsec-geodesic} the constant $0 < \e_L \leq \ez$ was defined
so that there is a projection map $\Pi_L \colon  \cN(L , \e_L) \to L$ along
the transverse geodesic $\e_L$-disks to $L$.

 Next, recall that  Proposition~\ref{prop-length} shows that
for  any leafwise $C^1$-path $\gamma$, we have the upper bound
$\ds    ||\gamma ||_{\cP} ~  \leq ~  \lceil (||\gamma ||_{\F}/\ez) \rceil$
for the number of plaques required to cover $\gamma$.

Suppose that $L$ is a compact leaf, and $x \in L$ is a fixed basepoint, then
for any $y \in L$ there is a leafwise geodesic $\gamma_{x,y} \colon [0,1]
\to L$ from $x$ to $y$ with  $\|\gamma \|_{\F} \leq diam(L)$. Thus,
$\gamma_{x,y}$ can be covered by at most $\lceil (||\gamma ||_{\F}/\ez)
\rceil$ plaques.

Recall that for $n > 0$, the number $ \mu(n, \ez) \leq  \ez$   was defined
in \eqref{def-mu}, and the constant $C_T  \geq 1$ was  introduced in
\eqref{eq-CT1}  and  \eqref{eq-CT2} as a bound on the   distortion of the
projection maps $\ds
\Pi_{\alpha_0}^{\F} \colon \cN(x,\ez) \to \cT_{\alpha_0}$.
Introduce the function
\begin{defn}\label{def-eps}
 For $0 < \e \leq \ez$ and $r > 0$,
\begin{equation} \label{def-Delta}
   \Delta(r, \e)  \equiv    \mu(\lceil r/\ez\rceil +2, \e/C_T)/C_T \ .
   \end{equation}
\end{defn}
 We scale both the domain variable $\e$ and the range value of $\mu$ by
$C_T$ so that we have uniform estimates for pairs of points in any geodesic
normal ball in the chart, a fact which will be used later.

\begin{lemma}\label{lem-eps}
Given $0 < \e \leq \ez$, and a leafwise $C^1$-path $\gamma \colon [0,1] \to
L$ of length at most $r$,   the transverse holonomy along $\gamma$ defines a
smooth embedding
$$\hh_{\gamma} \colon {\cN}(\gamma(0), \Delta(r,\e)) \to {\cN}(\gamma(1),
\e) ~, ~ \hh_{\gamma}(\gamma(0)) = \gamma(1) \ . $$
\end{lemma}
\proof
Let $\ds \cP = \{\cP_{\alpha_0}(z_0), \cP_{\alpha_1}(z_1), \ldots ,
\cP_{\alpha_n}(z_n) \}$ for  $n = ||\gamma ||_{\cP}$ be a covering of
$\gamma$ by plaque chains with $z_i \in  \cT_{\alpha_i}$ for $0 \leq i \leq
n$,  as given by  Proposition~\ref{prop-length}.
Set $x = \gamma(0)$ and $y = \gamma(1)$, then $x \in \cP_{\alpha_0}(z_0)$
and $y \in \cP_{\alpha_n}(z_n)$.

Recall that the constant $C_T$ was chosen so  that for  the projection  $\ds
\Pi_{\alpha_0}^{\F} \colon \cN(x,\ez) \to \cT_{\alpha_0}$, and for $x' \in
\cN(x,\ez)$,  the condition \eqref{eq-CT2} implies that
\begin{equation}\label{eq-estimatelocal1}
 {\bf d}_{M}(x , x')/C_T ~  \leq ~  {\bf
d}_{\cT}(\Pi_{\alpha_0}^{\F}(x),\Pi_{\alpha_0}^{\F}(x')) ~ \leq  ~ C_T\,
{\bf d}_{M}(x,x') \ .
\end{equation}
Likewise, for $y' \in   \cN(y,\ez)$ we have
\begin{equation}\label{eq-estimatelocal2}
{\bf d}_{M}(y , y')/C_T ~  \leq ~  {\bf
d}_{\cT}(\Pi_{\alpha_n}^{\F}(y),\Pi_{\alpha_n}^{\F}(y')) ~ \leq  ~ C_T\,
{\bf d}_{M}(y,y') \ .
\end{equation}

Given $x' \in \cN(x,\Delta(r,\e))$ then by \eqref{eq-estimatelocal1}, for
$x_0 = \Pi_{\alpha_0}^{\F}(x) \in \cT_{\alpha_0}$ and $x_0' =
\Pi_{\alpha_0}^{\F}(x') \in \cT_{\alpha_0}$, we have
${\bf d}_{\cT}(x_0, x_0') \leq \mu(\lceil r/\ez\rceil +2, \e/C_T)$.

Then by Lemma~\ref{lem-ext} and the inclusion \eqref{eq-unifn}, we have that
  $\ds d_{\cT}(\hh_{\cP}(x), \hh_{\cP}(x')) \leq \e/C_T$.

  Observe that $\hh_{\cP}(x) = \Pi_{\alpha_n}^{\F}(\gamma_1)$ and so we set
$\hh_{\gamma}(x') =  (\Pi_{\alpha_n}^{\F})^{-1}(\hh_{\cP}(x'))$.
  Then condition \eqref{eq-CT2} implies that   $d_M(\hh_{\gamma}(x),
\hh_{\gamma}(x')) \leq \e$, and so $\hh_{\gamma}(x') \in \cN(y,\e)$ as was
to be shown.
\endproof

We apply  Lemma~\ref{lem-eps} to obtain the following fundamental property
of compact foliations.
\begin{prop}\label{prop-cptcov}
Let $L'$ be a compact leaf of $\F$.
Given $\Lambda >0$, there exists $0 <  \delta_{\Lambda} < \e_{L'}$ so that
if
     $L$ is  a compact  leaf with volume ${\rm vol}(L) < \Lambda$ and
$L \cap \cN(L',\delta_{\Lambda}) \not= \emptyset$, then
    $L \subset \cN(L' ,\e_{L'})$.
  \end{prop}
\proof

Let $R(\Lambda)$ be the constant introduced in the proof of
Proposition~\ref{prop-volume}, so that
  ${\rm vol}(L) < \Lambda$  implies that  $\diam (L) \leq R(\Lambda)$.
   Set   $\delta_{\Lambda}  = \Delta(R(\Lambda), \e_{L'})/2$.
 Assume that $L \cap \cN(L',\delta_{\Lambda}) \not= \emptyset$, then  there
exists   $x' \in L'$ such that there exists
    $ x \in L \cap \cN(x',\delta_{\Lambda})$. We must  show that  $L
\subset  \cN(L',\e_{L'})$.

We have  $ x \in L \cap \cN(x',\delta_{\Lambda})$ and so $D_M(x',x) <
\delta_{\Lambda}$. Then  by assumption (4) in Section~\ref{subsec-geodesic},
 there exists $w' \in L' \cap \cN(x, 2 \delta_{\Lambda})$, and so
$d_M(x,w') < 2\delta_{\Lambda}  = \Delta(R(\Lambda), \e_{L'})$.

Let $z \in L$, then there is a leafwise geodesic path $\gamma_{x,z} \colon
[0,1] \to L$ with $\|\gamma_{x,Z}\|_{\F} \leq R(\Lambda)$ and
$x=\gamma_{x,z}(0)$ and $z = \gamma_{x,z}(1)$.
By Lemma~\ref{lem-eps}, there exists a holonomy map $\ds \hh_{\gamma_{x,z}}
\colon {\cN}(x, 2\delta_{\Lambda}) \to {\cN}(z, \e/2)$.

As $w' \in L' \cap \cN(x, 2 \delta_{\Lambda})$ so $w'' =
\hh_{\gamma_{x,z}}(w')$ is well-defined. If $\cP$ is the plaque chain
containing $\gamma_{x,z}$ chosen in the proof of Lemma~\ref{lem-eps},
then there is a  corresponding plaque chain $\ds \cP' =
\{\cP_{\alpha_0}(z_0'), \cP_{\alpha_1}(z_1'), \ldots , \cP_{\alpha_n}(z_n')
\}$ in $L'$ from $w'$ to $w''$, where each $z_i' \in \cT_{\alpha_i}$. We say
that $\cP'$ \emph{shadows} in $L'$ the plaque chain $\cP$ in $L$.
In particular, we have $w'' \in \cP_{\alpha_n}(z_n') \subset L'$.

 The conclusion of Lemma~\ref{lem-eps} is that for $2\delta_{\Lambda}  =
\Delta(R(\Lambda), \e_{L'})$ then $w'' \in \cN(z, \e_{L'})$, and hence
$d_M(z,w'') < \e_{L'}$. Thus, $z \in \cN(L', \e_{L'})$ for all $z \in L$,
hence  $L \subset \cN(L' ,\e_{L'})$.
\endproof

Proposition~\ref{prop-cptcov} has the following useful consequence.

\begin{cor}\label{cor-cptcov} Let $L_0$ be a compact leaf of $\F$.
Given $\Lambda >0$, there exists $0 <  \delta_{\Lambda}  < \e_{L_0}$ so that
if
     $L_1$ is  a compact  leaf with volume ${\rm vol}(L_1) < \Lambda$ and
    $L_1 \cap \cN(L_0,\delta_{\Lambda}) \not= \emptyset$, then $L_1 \subset
\cN(L_0, \e_{L_0})$.
    Moreover, the   projection map
    $\Pi_{L_0} \colon \cN(L_0, \e_{L_0}) \to L_0$ restricted to $L_1$ is a
covering map onto $L_0$.  Furthermore, if the tangent bundle $T\F$ is
orientable, then
    $vol (L_1) \leq 2 \, d_* \, vol(L_0)$
    where $d_*$ is the homological degree of the covering map $\Pi_{L_0}
\colon L_1 \to L_0$.
  \end{cor}
\proof
Let $\delta_{\Lambda}$ be as defined in Proposition~\ref{prop-cptcov}, then
$L_1 \subset \cN(L_0, \e_{L_0})$ follows.

By the assumption (4) in Section~\ref{subsec-geodesic}, the leaves of $\F$
are uniformly transverse to the fibers of $\Pi_{L_0} \colon \cN(L_0,
\e_{L_0}) \to L_0$, so the restriction to $L_1$ is a covering map. As $L_0$
and $L_1$ are compact, the map $\Pi_{L_0} \colon  L_1 \to L_0$ is onto.
Assume that the tangent bundle $T\F$ is oriented, then we can choose a
positively-oriented Riemannian volume form on the leaves of $\F$, whose
restriction to a leaf $L$ is denoted by $\omega_L$. Then the closed $p$-form
$vol(L)^{-1} \cdot
\omega_L$ is   dual to  the fundamental class $[L]$ of $L$.

The homological degree $d_*$ of a covering map equals its covering degree,
and is given by
\begin{equation}\label{eq-degrees}
d_*    = d_* \cdot  vol(L_1)^{-1} \cdot \int_{L_1} \ \omega_{L_1}  =
vol(L_0)^{-1} \cdot \int_{L_1} \ \Pi_{L_0}^*(\omega_{L_0}) \ .
\end{equation}

Condition (6) of Section~\ref{subsec-geodesic} gives that for $x' \in L_1$
and $x = \Pi_{L_0}(x') \in L_0$ we have
\begin{equation}\label{eq-pointwisebounds}
1/2 \cdot   \omega_{L_1}|_{x'}  ~ \leq ~  \Pi_{L_0}^*(\omega_{L_0})|_x   ~
\leq ~ 2 \cdot   \omega_{L_1}|_{x'} \ ,
\end{equation}
and thus combining  \eqref{eq-degrees} and \eqref{eq-pointwisebounds} we
obtain
\begin{equation}
 1/2 \cdot vol(L_1) = 1/2 \cdot \int_{L_1}\omega_{L_1} \leq d_* \cdot
vol(L_0) = \int_{L_1} \ \Pi_{L_0}^*(\omega_{L_0}) \leq 2 \cdot \int_{L_1} \
\omega_{L_1} = 2 \cdot vol(L_1) \ .
\end{equation}
Thus, $\ds vol(L_1)  \leq 2 d_* \cdot vol(L_0)$, as was to be shown.
\endproof

\section{Properties of compact foliations}\label{sec-properties}

In this section, $\F$ is assumed to be a compact foliation of a   manifold
$M$ without boundary. The geometry of compact foliations has been studied
by Epstein \cite{Epstein1972, Epstein1976}, Millett \cite{Millett1974}, Vogt
\cite{Vogt1976, Vogt1977a, Vogt1977b} and Edwards, Millett and Sullivan
\cite{EMS1977}. We recall some of their results.

\subsection{The good and the bad sets} \label{subsec-goodbad}

 Let $vol(L)$ denote the volume of a leaf $L$ with respect to the
Riemannian metric induced from $M$.  Define the volume function on $M$ by
setting   $v(x) = vol(L_x)$. Clearly, the function $v(x)$ is constant along
leaves of $\F$, but   need not be continuous on $M$.

The {\it bad set} $X_1$ of $\F$ consists of the points $y \in M$ where the
function $v(x)$ is not   bounded in any open neighborhood of $y$.
By its definition, the bad set $X_1$ is saturated. Note also that
$$X_1 = \bigcup_{n=1}^{\infty} \;  X_1 \cap vol^{-1}(0,n] \ .$$
The leaves in the intersection  $X_1 \cap vol^{-1}(0,n]$  have  volume at
most  $n$, while  $v(x)$ is not locally bounded in any open neighborhood of
$y \in X_1$, therefore each set   $X_1 \cap vol^{-1}(0,n]$ has no interior.
By the Baire category theorem,  $X_1$ has no interior.

The  complement
$G = M \smallsetminus  X_1$ is called the {\it good set}. The   holonomy of
every leaf
$L \subset G$ is finite,  thus by the Reeb Stability Theorem, $L$ has an
open saturated neighborhood consisting of leaves with finite holonomy.
Hence,  $G$ is   an open set,  $X_1$ is closed, and   the   leaf space
$G/\F$ is Hausdorff.

 Inside the good set is  the open dense saturated subset  $G_e \subset G$
consisting  of leaves without holonomy. Its complement $G_h = G
\smallsetminus G_e$
consists of leaves with non-trivial finite holonomy.

\subsection{The Epstein filtration}\label{subsec-epstein}

The restriction of the volume function $v(x)$ to $X_1$ again need not be
locally bounded, and the construction of  the bad set can be iterated to
obtain the  {\it Epstein filtration}:
$$M = X_0 \supset X_1 \supset X_2 \supset \cdots \supset X_{\alpha} \supset
\cdots  \ . $$
The definition of the sets $X_{\alpha}$ proceeds inductively: Let $   \alpha
> 1$ be a countable  ordinal, and assume that $X_{\beta}$ has been defined
for $\beta < \alpha$.  If $\alpha$ is a successor ordinal,  let $\alpha =
\gamma + 1$ and  define $X_{\alpha}$ to be the closed saturated set of   $y
\in X_{\gamma}$ where the function $v(x)$ is not   bounded in any relatively
open neighborhood  of $y \in X_{\gamma}$ in $X_{\gamma}$.

If $\alpha$ is a limit ordinal, then define
$\ds  X_{\alpha} = \bigcap_{\beta < \alpha} X_{\beta}$.

For $\beta < \alpha$, the set  $X_{\alpha}$ is   nowhere dense   in
$X_{\beta}$. Note that since each set $M \smallsetminus X_{\alpha}$ is open,
the
filtration is at most countable.   The {\it filtration length} of $\F$ is
the ordinal $\alpha$ such that  $X_{\alpha} \not= \emptyset$ and $X_{\alpha
+ 1}  = \emptyset$.

  Vogt \cite{Vogt1977b}  showed  that for any finite ordinal $\alpha$, there
is a compact foliation of a compact manifold with filtration length
$\alpha$.   He also  remarked that given any countable ordinal $\alpha$,
the construction can be modified to produce a foliation with filtration
length $\alpha$.
Such examples show that the bad set $X_1$ and the subspaces $X_{\alpha}$
need not be  finite unions or intersections of submanifolds;  they may have
pathological topological structure, especially when the filtration length is
an infinite ordinal.

\subsection{Regular  points}\label{subsec-regular}

 A point   $x \in X_1$ is called a \emph{regular point}    if the
restricted  holonomy of $\F|X_1$ is trivial at $x$.  Equivalently, the
regular points are the points of continuity for the restricted volume
function $v|X_1\colon X_1  \to \mR$.
If $X_1 \not= \emptyset$, then the regular  points form an open and dense
subset of $X_1 \smallsetminus X_2$.
We recall a key result of Edwards, Millett, and Sullivan (see \S 5 of
\cite{EMS1977}.)
\begin{prop}[Moving Leaf]\label{prop-mlp}
Let $\F$ be a compact foliation of an oriented manifold $M$ with orientable
normal bundle.
Suppose that $X_1$ is compact and non-empty. Let $x \in X_1$ be a regular
point. Then there exists a   leaf $L \subset G_e$  and a smooth
isotopy $h \colon L \times [0,1) \to G_e$ such that:
\begin{itemize}
\item For all $0 \leq t < 1$, $h_t \colon L \to L_t \subset M$ is a
diffeomorphism onto its image $L_t$
\item   $L_x$ is in the closure of the leaves ~~ $\ds  \bigcup_{t>1- \delta}
L_t$ ~~ for any $\delta > 0$
\item   $\ds  \limsup_{t\to 1} ~ vol(L_t) = \infty$
\end{itemize}
\end{prop}

While the ``moving leaf'' $L_t$   limits on    $X_1$,    the moving leaf
cannot accumulate on a single compact leaf  of $X_1$. This follows because a
compact leaf $L$ admits a relative homology dual cycle, which for $\e > 0$
sufficiently small and $x \in L$, is represented by the transverse disk
${\bf B}_{\cT}(x, \e)$. This disk  intersects $L$ precisely in the point
$x$, hence
 the relative homology class $[{\bf B}_{\cT}(x, \e), \partial {\bf
B}_{\cT}(x, \e)]$ is Poincar\'e dual to the fundamental class $[L]$.
Assuming that
  $\{L_t\}$   limits on    $L$, for $t < 1$ sufficiently close to $1$, each
$L_t \subset {\cN}(L,\e)$ and so  the intersection number $[L_t \cap {\bf
B}_{\cT}(x, \e)] = [L_t] \cap [{\bf B}_{\cT}(x, \e), \partial {\bf
B}_{\cT}(x, \e)]$ is constant.  Thus   the   leaves $\{L_t\}$ have   bounded
volume as $t \to 1$,  which is a contradiction.

It is precisely this ``non-localized limit behavior'' for leaves with
unbounded volumes approaching the bad set   which makes the study of compact
foliations with non-empty bad sets so interesting, and difficult. There seem
to be  no   results in the literature describing how these paths of leaves
must behave in the limit.

\subsection{Structure of the good set}\label{subsec-good}
  Epstein \cite{Epstein1976} and Millett \cite{Millett1974} showed that for
a  compact foliation $\F$ of a manifold $V$, then
  $$ v(x) ~ {\rm is} \, {\rm locally} \, {\rm bounded}  \Leftrightarrow V/\F
~ {\rm is }\, {\rm   Hausdorff} \Leftrightarrow
 ~ {\rm  the} \, {\rm  holonomy} \, {\rm of}\,  {\rm every} \, {\rm  leaf}
\, {\rm  is} \, {\rm finite} \ . $$

By definition, the  leaf volume function is locally bounded on  the good set
$G$, hence the restriction of $\F$ to $G$ is   compact Hausdorff, and all
leaves of $\F |  G$ have finite holonomy group. Epstein and Millett showed
there is a much more precise structure theorem for the foliation $\F$ in an
open neighborhood of a leaf of the good set:

 \begin{prop}\label{prop-locstruct}  Let $V$ denote an open connected
component of the good set $G$, and $V_e = V \cap G_e$ the set of leaves with
no holonomy.
 There exists a    ``generic leaf'' $L_0 \subset V_e$,  such that  for each
$x \in V$ with leaf $L_x$   containing $x$,
\begin{enumerate}
\item   there is a   finite subgroup $H_x$ of the orthogonal group $\bf
O(q)$
and a free action  $\alpha_x$ of $H_x$  on $L_0$
\item  there exists a diffeomorphism of the twisted product
\begin{equation}\label{eqn.rep}
  \phi_x \colon L_0  \times_{H_x} \D \to V_x
\end{equation}
onto  an open saturated neighborhood $V_x$ of $L_x$ (where $\D$ denotes the
unit disk in ${\mR}^q$)
\item the diffeomorphism  $\phi_x$ is leaf preserving, where $\ds    L_0
\times_{H_x} \D$~ is foliated by the images  of $L_0 \times \{w\}$ for ~ $w
\in \D$ under the quotient map
$\ds  {\cQ} \colon L_0  \times  \D \to L_0  \times_{H_x} \D$
\item   $\phi_x$   maps   $L_0/H_x \cong  L_0 \times_{H_x} \{0\}$
diffeomorphically to $L_x$
\end{enumerate}
In particular, if $x \in V_e$ then $H_x$ is trivial, and $\phi_x$ is a
product structure for a neighborhood of $L_x$.
\end{prop}

The open set $V_x$ is called a standard neighborhood of $L_x$, and the
4--tuple  $(V_x,\phi_x,H_x,\alpha_x)$ is called a {\it standard local model}
for $\F$. Note that, by definition,     $V_x \subset G$  hence $V_x \cap X_1
= \emptyset$.

The  Hausdorff space  $G/{\F}$ is a Satake manifold;  that is,
for each point $b \in G/\F$ and   $\pi(x) = b$
the leaf $L_x$ has an  open foliated neighborhood $V_x$ as above,
and
$\ds   \phi_x \colon L_0  \times_{H_x} {\D} \to V_x $
induces a coordinate map     $\ds  \varphi_b \colon {\D}/H_x \to W_b$, where
$W_b = \pi(V_x)$.
The open sets  $W_b \subset G/{\F}$ are called  basic  open sets for $G/\F$.
Note  also that   $\pi$ is a closed map \cite{Epstein1976,Millett1974}.

\section{Properties of foliated homotopies}\label{sec-homotopy}

 In this section, we study some of the geometric and topological properties
of a foliated homotopy of a compact leaf.
 These results play an essential role in our proof of
Proposition~\ref{prop-npp},      and  hence of   Theorem~\ref{thm-compact}.
The main result of this section is the following,  which yields   an upper
bound on both the volumes of the compact leaves,   and
 the topological degrees of the covering maps,  which arise in a homotopy of
a compact leaf.  Note that the results of this section apply to any
$C^1$-foliation of a manifold $M$.
  \begin{prop}\label{prop-bddvol3}
 Let $\F$ be a $C^1$ foliation of a   manifold $M$,
let  $L$ be a compact leaf, and let $H \colon L \times [0,1] \to M$ be a
foliated homotopy, and let $L_t$ denote the compact leaf containing
$H_t(L)$.
Assume that   both the tangent  bundle $T\F$ and the normal bundle $Q$ to
$\F$ are oriented. Then
 there exists  $d_* > 0$, depending on $H$ and $L$,  such that
\begin{equation} \label{eq-unideg}
1 \leq \deg(H_t \colon L_0 \to L_t) \leq d_* \ .
\end{equation}
Moreover, there exists an integer
$k \geq 0$  such that
\begin{equation} \label{eq-unifbdd}
vol(L_t) \leq 4^k d_* \cdot vol(L)  ~ , \quad {\rm   for ~ all} ~ 0 \leq t
\leq 1 \ .
\end{equation}
\end{prop}
\proof

We first recall a result which     holds for all $C^1$-foliations.
\begin{thm}\cite[Corollary 1.4]{Hurder2006}\label{thm-compactleaf}
 Let $\F$ be a $C^1$ foliation of a compact manifold $M$. Let  $L$ be a
compact leaf, and $H \colon L \times [0,1] \to M$ be a foliated homotopy.
Then for all $0 < t \leq 1$, the image $H_t(L)$ is contained in a compact
leaf $L_t$ of $\F$, and moreover, the map $H_t \colon L \to L_t$ is
surjective.
\end{thm}

Now let $0 \leq t \leq 1$, and let $L_t$ be the compact leaf for which
$H_t(L) \subset L_t$. We set $L_0 = L$ for notational consistency. Then
there exists $0< \e_t = \e_{L_t} \leq \ez$ such that we have a normal
$\e_t$-bundle projection map
$\Pi_{L_t}  \colon \cN(L_t,\e_t) \to L_t$. The subset $\cN(L_t,\e_t) \subset
M$ is open, and $H$ is uniformly continuous, so there exists $\delta_t > 0$
such that $H_s(L) \subset \cN(L_t,\e_t)$ for all $t - \delta_t < s < t +
\delta_t$. For such $s$, the map $H_s \colon L_0 \to L_s$ is onto, so the
leaf $L_s \subset \cN(L_t,\e_t)$, hence the restriction $\Pi_{L_t} \colon
L_s \to L_t$ is a covering map.

The maps $H_s, H_t \colon L_0 \to \cN(L_t,\e_t)$ are homotopic in
$\cN(L_t,\e_t)$, hence for their induced maps on fundamental classes their
degrees satisfy
\begin{equation}
\deg(H_t \colon L_0 \to L_t)  = \deg(\Pi_{L_t} \colon L_s \to L_t)  \cdot
\deg(H_s \colon L_0 \to L_s) \ .
\end{equation}
The homological degree of a covering map equals its covering degree, thus
the covering degree of $\Pi_{L_t} \colon L_s \to L_t$ divides the
homological degree of $\deg(H_t \colon L_0 \to L_t)$.

The collection of open intervals $\{\cI_t = (t-\delta_t, t+\delta_t) \mid 0
\leq t \leq 1\}$ is an open covering of $[0,1]$, so there exists a finite
set $\{0 = t_0 < t_1 < \cdots < t_{k-1} < t_k =1\}$ so that the collection
$\{\cI_{t_i} \mid 0 \leq i \leq k\}$ is a finite covering of $[0,1]$. Choose
a sequence $\{0 < s_1 < \cdots < s_{k-1} < 1\}$ such that
$$ t_{\ell -1} < s_{\ell} < t_{\ell} \quad , \quad t_{\ell} -
\delta_{t_{\ell}} < s_{\ell} < t_{\ell -1} + \delta_{t_{\ell-1}}$$
and hence  $s_{\ell} \in \cI_{t_{\ell -1}} \cap \cI_{t_{\ell}}$. Thus for
the choices of the constants $\delta_t$ for each $0< \ell < k$, we have the
inclusions
$\ds L_{s_{\ell}} \subset \cN(L_{t_{\ell -1}},\e_{t_{\ell -1}}) \cap
\cN(L_{t_{\ell}},\e_{t_{\ell}})$. Thus, there are finite covering maps
\begin{equation}\label{eq-commensurable}
 \Pi_{L_{t_{\ell -1}}} \colon L_{s_{\ell}} \to L_{t_{\ell -1}} ~ , ~
\Pi_{L_{t_{\ell}}} \colon L_{s_{\ell}} \to L_{t_{\ell}} \quad , ~ {\rm for ~
each} ~ 1 \leq \ell < k-1 \ .
\end{equation}
The collection of maps  \eqref{eq-commensurable} is called a \emph{geometric
correspondence} from $L_0$ to $L_1$. We have shown:

\begin{lemma}\label{lem-commen}
Let $\F$ be a $C^1$ foliation of  $M$, $L$ a compact leaf of $\F$, and
 $H \colon L \times [0,1] \to M$   a foliated homotopy.
 Then there exists a geometric correspondence from $L_0 = L$ to $L_1 =
H_1(L)$.
\end{lemma}

Introduce the following integer constants associated to a correspondence
\eqref{eq-commensurable}, for $0< \ell < k$:
\begin{eqnarray}
a_{\ell} & = &  \deg(\Pi_{L_{t_{\ell -1}}} \colon L_{s_{\ell}} \to
L_{t_{\ell -1}} )  \label{eq-correpondence-aell}\\
b_{\ell} & = &  \deg(\Pi_{L_{t_{\ell}}} \colon L_{s_{\ell}} \to
L_{t_{\ell}})   \label{eq-correpondence-bell}\ .
\end{eqnarray}
Note that $a_{\ell}$ and $b_{\ell}$ are equal to the covering degrees of the
covering maps in \eqref{eq-correpondence-aell} and
\eqref{eq-correpondence-bell}, and that $a_1 =1$ as the leaf $L_{s_1}$ must
be a diffeomorphic covering of $L_0$.
Then  the choice of each $\e_t \leq \ez$ we can apply the estimate
\eqref{eq-univol} as in the proof of Corollary~\ref{cor-cptcov} to obtain,
for $1 \leq \ell < k$,
$$
\begin{array}{lcccl}
a_{\ell} /2 \cdot vol(L_{t_{\ell -1}}) & \leq & vol (L_{s_{\ell}}) & \leq &
2   a_{\ell} \cdot vol(L_{t_{\ell -1}}) \\
b_{\ell} /2 \cdot vol(L_{t_{\ell}}) & \leq & vol (L_{s_{\ell}}) & \leq & 2
b_{\ell} \cdot vol(L_{t_{\ell}})
\end{array}
$$
We can then combine these sequences of upper and lower estimates to obtain
the estimate:
 \begin{equation} \label{eq-volest3}
4^{-k} \, \frac{b_1 \cdots b_{k-1}}{a_1 \cdots a_{k-1}}  \cdot vol(L_0) ~
\leq ~  vol (L_{1}) ~ \leq  ~ 4^k \, \frac{ a_1 \cdots a_{k-1}}{ b_1 \cdots
b_{k-1}}   \cdot vol(L_{0}) \ .
 \end{equation}
Set $\ds d_* =  \frac{a_1 \ldots a_{k-1}}{b_1 \ldots b_{k-1}}$ and we obtain
the estimate \eqref{eq-unifbdd} for $t =1$.
The uniform bound \eqref{eq-unideg} follows from the argument above,
considering only the homological degrees of the covering maps and ignoring
the volume estimates. A minor modification of the above arguments also
yields these estimates for the   values   $0 < t < 1$.
This completes the proof of Proposition~\ref{prop-bddvol3}.
\endproof

\section{Tame points in the bad set}\label{sec-tame}

In this section, we introduce  the concept of a ``tame point'' in the bad
set $X_1$, which is a point $x \in X_1$ that can be approached by a path in
the good set. The main result of this section  proves  the existence of tame
points, using a more careful analysis of the ideas of the Moving Leaf
Proposition~\ref{prop-mlp}. Tame points are used in  section 6 for  studying
the deformations of the bad set under foliated homotopy.

 Recall that the  bad set  $X_1$  consists of the points $y \in M$ where the
leaf volume function $ v(x)$ is not   bounded in any open neighborhood of
$y$, and  is closed, saturated and has no interior. A point  $x_1 \in X_1$
is \emph{regular} if   the  restriction of the leaf volume function $v
\colon X_1 \to {\mathbb R}^+$   is   continuous at $x_1$. Equivalently, $x_1
\in X_1$ is a regular point if  the   holonomy of the restriction of $\F$ to
$X_1$ is trivial in some relatively open neighborhood of $x_1 \in X_1$.

   \begin{defn}\label{def-pp}
      A regular point $x_1 \in X_1$ is   {\em tame} if there exists $\e > 0$
and a   transverse  $C^1$-path
  \begin{equation}\label{eq-tame}
   \gamma \colon [0,1] \to  \left(\cN(x_1,\e) \cap G_e\right) \cup \{x_1\}
  \end{equation}
with  $\gamma(t)  \in G_e$ for $0\leq t < 1$, $\gamma(1) = x_1$ and such
that $v(\gamma(t)))$ tends uniformly to infinity as $t \to 1$.
   \end{defn}

  Let $X_1^t \subset X_1$ denote the subset of tame points.

   Since the restricted path     $\gamma \colon [0,1) \to    G_e$ lies in
the set of leaves without holonomy,
   it follows that  there is a   foliated isotopy $\Gamma  \colon
L_{\gamma(0)} \times [0,1)  \to    G_e$ such that $\Gamma_t(\gamma(0)) =
\gamma(t)$.
 Thus,   a tame point $x$ is directly approachable by a family of moving
leaves whose volumes tend uniformly  to infinity.

In the examples constructed by Sullivan   \cite{Sullivan1976},  it is easy
to see that every regular  point  is a tame point.
 In general, though,
      Edwards, Millet, and Sullivan   specifically point out that their
proof  of the Moving Leaf Proposition~\ref{prop-mlp}
      in \cite{EMS1977}  does not claim that a regular point is a tame
point. The problem is due to the possibility that the complement of the bad
set need not be locally connected in a neighborhood of a point in the bad
set.
        In their proof, the moving leaf is defined by a curve that follows
``an end of the good set'' out to infinity, passing through points where the
volume is tending to infinity along the way. This end of the good set is
contained in arbitrarily small $\e$-neighborhoods of the bad set, but they
do not control the behavior of the end. Thus, the existence of a tame point
is  asserting the existence of a ``tame end'' of the good set on which the
volume function is unbounded, and which is defined by open neighborhoods of
some point in the bad set.

 \begin{prop}\label{prop-tame}
 Let $\F$ be a compact, $C^1$-foliation of a manifold $M$, and assume that
the tangent  bundle $T\F$ and the normal bundle $Q$ to $\F$ are oriented.
Then the set of tame points  $X_1^t$ is dense in  $X_1$.
 \end{prop}
\proof
The proof of this result involves several technical steps, so we first give
an overview of the strategy of the proof. Let $x_1 \in X_1$ be a regular
point, and $L_{1}$ the leaf through $x_1$. We use a key result in the proof
of the Moving Leaf Lemma to obtain an open neighborhood $U$ of $x_1$ in its
transversal space, on which the volume function is unbounded. We then choose
a regular point $x_* \in U \cap X_1$ which is sufficiently close to $x_1$,
so that the leaf $L_*$ through $x_*$ is a diffeomorphic covering of the leaf
$L_1$. Moreover, the point $x_*$ is approachable by a path in the good set.
Then we argue by contradiction, that if the leaf volume function does not
tend uniformly to infinity along this path, then each leaf through a point
in  the   set $U \cap G$ is also a covering space of $L_1$ with uniformly
bounded covering degree, from which we conclude that the volume function is
bounded on the leaves through points in $U$, contrary to   choice. It
follows that $x_*$  is a tame point   which is arbitrarily close to $x_1$.
The precise proofs of these claims requires
that we first  establish some technical properties of the foliation $\F$
in a normal neighborhood of $L_1$.

The leaf $L_1$ is compact,  hence has finitely-generated fundamental group.
Thus, we can choose  a finite generating set
 $\{[\tau_1] , \ldots , [\tau_k]\}$ for
 $\pi_1(L_1, x_1)$, where $[\tau_i]$ is represented by a smooth closed path
$\tau_i \colon [0,1] \to L_1$
    with basepoint $x_1$. Let $ \| \tau_i \|$ denote the path length of
$\tau_i$.
Then set
\begin{equation} \label{eq-frakp}
D_{L_1} = 2 \max \left\{  \diam (L_1),  \| \tau_1 \|, \ldots , \|\tau_k\|
\right\} \ .
\end{equation}

 Recall that in Section~\ref{subsec-geodesic}, given a compact leaf $L$  the
constant $0 < \e_L \leq \ez$ was defined  so that there is a projection map
$\Pi_L \colon  \cN(L , \e_L) \to L$ along the transverse geodesic
$\e_L$-disks to $L$. Set $\eone = \e_{L_1}$  so that the normal projection
map    $\Pi_{L_1}  \colon \cN(L_1, \eone) \to L_1$ is well-defined.
Then  set $\e_2 =  \Delta( D_{L_1}, \eone)$ where $ \Delta(D_{L_1}, \eone)$
is
defined in Definition~\ref{def-eps}.
Then  by  Lemma~\ref{lem-eps},  for any path $\sigma \colon [0,1] \to L_1$
with $\sigma(0) = x_1$ and $\| \sigma \| \leq D_{L_1}$ the transverse
holonomy
maps are defined for all $0 \leq t \leq 1$,
\begin{equation}\label{eq-lift2}
\hh_{\sigma} \colon \cN(x_1, \etwo) \to   \cN(\sigma(t), \eone) \ .
\end{equation}
In particular,
the  holonomy map $\hh_i$   along each closed path $\tau_i$   is defined on
the transverse disk $\cN(x_1, \etwo)$.
That is, the transverse holonomy along $\tau_i$ is represented by a local
homeomorphism  into
\begin{equation}\label{eq-lift1}
\hh_i \colon \cN(x_1, \etwo) \to   \cN(x_1, \eone) \ .
\end{equation}

The assumption that $x_1 \in X_1$ is a regular point implies that  the
germinal  holonomy at $x_1$ of the restricted foliation $\F|X_1$ is trivial.
 Thus we can choose      $0 < 2 \delta \leq \etwo$   sufficiently small,
so that each holonomy map
$\hh_i$ restricted to    $X_1 \cap \cN(x_1, 2\delta)$ is the identity map.
It follows that  the    holonomy of $\F$ restricted to  the closure
\begin{equation}\label{eq-choosedelta}
Z_1 = \overline{X_1 \cap \cN(x_1, \delta)} = X_1 \cap \overline{\cN(x_1,
\delta)}  \subset  X_1 \cap \cN(x_1, 2\delta)
\end{equation}
is trivial. Hence,   every point in $Z_1$ is a regular point of the bad set.
It follows that the saturation $Z_{\F}$ of $Z_1$ is   a fibration over the
closed set $Z_1$, and that
the  leaf volume function $v(y)$  is uniformly continuous and hence bounded
on the compact set $Z_1$.
Thus, we may   assume that $\delta $ is sufficiently small  so  that
$Z_{\F} \subset \cN(L_1, \eone)$.
That is, for each $z \in Z_1$ the leaf $L_z \subset \cN(L_1, \eone)$.

Next consider the properties of the  the normal projection $\Pi_{L_1} \colon
\cN(L_1, \eone) \to L_1$ when restricted to leaves in $\cN(L_1, \eone)$.
 The restriction  $\pi^z \equiv \Pi|L_z \colon L_z \to L_1$      is a
covering projection, which is a diffeomorphism as $\F | Z_1$ has no
holonomy, and
by  the assumption that $\eone \leq \ez$, the estimate \eqref{eq-unimetric}
implies the map $\pi^z$ is   a quasi-isometry with expansion constant
bounded by $2$.

Note that  $\cN(x_1, 2\delta)$ is contained in the the normal transversal
$\cN(x_1, \eone)$, so
by  definition for  $z \in Z_1$ we have $\pi^z(z) = x_1$, and thus given   a
path $\sigma \colon [0,1] \to L_1$  with $\sigma(0) = x_1$, there is a lift
 $\sigma^z \colon [0,1] \to L_z$  with $\sigma^z(0) = z$ and   $\pi^z \circ
\sigma^z(t) = \sigma(t)$ for all $0 \leq t \leq 1$.
In particular,   the closed loop $\tau_i$ lifts via $\pi^z$  to a closed
loop
$\tau_i^z \colon [0,1] \to L_z$ with endpoints $z$.
The homotopy classes of the lifts, $\{ [\tau_1^z] , \ldots , [\tau_k^z]\}$,
yield  a generating set for $\pi_1(L_z, z)$,
which   have  a uniform bound $\| \tau_i^z \| \leq  D_{L_1}$ on their path
lengths.

For an arbitrary point $y_0 \in \cN(x_1, \delta)$ and path    $\sigma \colon
[0,1] \to L_1$ with $\sigma(0) = x_1$ and path length $\| \sigma \| \leq
D_{L_1}$, by the choice of $\delta$  the   transverse holonomy map in
\eqref{eq-lift2} is defined at $y_0$ hence   there is a lift   of the path
$\sigma$ to a path
$\ds \sigma^y \colon [0,1] \to L_y \cap  \cN(L_1 , \eone)$
 with $\sigma^y(0) = y_0$ and  $\pi^y \circ \sigma^y(t) = \sigma(t)$ for all
$0 \leq t \leq 1$.
This lifting property need not  hold for paths longer than $D_{L_1}$, as
there
may be leaves of $\F$ which intersect  the normal neighborhood $\cN(L_1,
\delta)$ but are not contained in $\cN(L_1, \eone)$.

We observe  a technical point about the distances in the submanifold
$\cN(x_1, \eone) \subset M$. The inclusion $\cN(x_1, \eone) \subset M$
induces a Riemannian metric on $\cN(x_1, \eone)$ which then defines a
path-length distance function on this subspace. Unless $\cN(x_1, \eone)$ is
a totally geodesic submanifold of $M$, the induced distance function on
$\cN(x_1, \eone)$ need not agree with the restricted path-length  metric
from $M$.
  For $y \in \cN(x_1, \eone)$ and  $0 < \lambda \leq \eone$, let    $B_T(y,
\lambda) \subset \cN(x_1, \eone)$ denote the open ball of radius $\lambda$
about $y$  for the induced Riemannian metric on  $\cN(x_1, \eone)$.

Now consider an arbitrary point  $y \in \cN(x_1, \delta)$ and assume that
$L_y \subset \cN(L_1, \delta)$,    so that $L_y$ is in the domain of the
projection $\Pi_{L_1} \colon \cN(L_1, \eone) \to L_1$.
Given a path $\sigma \colon [0,1] \to L_y$ with $\sigma(0) = y$, then
Lemma~\ref{lem-eps}
implies that there exists $0 < \lambda < \delta$, which depends on the
length $\|\sigma\|$, so that
for $y' \in B_T(y, \lambda)$   there is a  path  $\sigma^{y'} \colon [0,1]
\to L_{y'}$  satisfying $\Pi_{L_y}(\sigma(t)) = \Pi_{L_y}(\sigma^{y'}(t))$
for $0 \leq t \leq 1$.
We call the path $\sigma^{y'}$ a lifting of $\sigma^{y}$ from $L_y$ to
$L_{y'}$.

   After these technical preliminary results, we begin the proof of
Proposition~\ref{prop-tame}.
 First, recall   a key fact   from the proof of the Moving Leaf
Proposition~\ref{prop-mlp}, whose proof was in turn based on ideas of
Montgomery \cite{Montgomery1937} and Newman \cite{Newman1931}.   (In
particular,   Figure~3  on page 23 of   \cite{EMS1977} and the arguments
following it are pertinent to the arguments below.)

\begin{lemma}\label{lem-unbounded}
For   $\delta > 0$ sufficiently small,  there is an open connected component
$U$ of $\cN(x_1,\delta) \smallsetminus Z_1$ on which the volume function
$v(y)$ is
unbounded on  the open neighborhood   $U \cap \cN(x_1 , \delta/2)$.
\end{lemma}

 Next, fix a choice of regular point $x_1 \in X_1$ and sufficiently small
constant $\delta > 0$ as above so that \eqref{eq-choosedelta} holds,  then
choose a point $y_0 \in U \cap \cN(x_1 , \delta/2)$.
Let $x_* \in Z_1$ be a closest point to $y_0$ for the induced metric on
$\cN(x_1,\delta)$.
That is, consider the sequence of closed balls $\ds \overline{B_{T}(y_0,
\lambda)} \subset \cN(x_1,\delta) \smallsetminus Z_1$   for $\lambda > 0$,
expanding
until there is a first contact with the  frontier of $U$, then  $x_*$ is
contained in this intersection.
Let $\delta_0 \leq \delta/2$ denote the radius of first contact, hence
$\delta_0$ equals the distance from $y_0$ to $x_*$ in the induced
path-length metric on $\cN(x_1, \eone)$. Then $B_T(y_0, \delta_0) \subset U$
and $x_* \in \overline{B_T(y_0, \delta_0)} \cap  Z_1$.  (This is illustrated
in Figure~\ref{fig:tamepoint} below.)  Let $L_* = L_{x_*}$ denote the leaf
containing $x_*$.

 \begin{figure}[htbp]
\begin{center}
\includegraphics[width=0.5\textwidth]{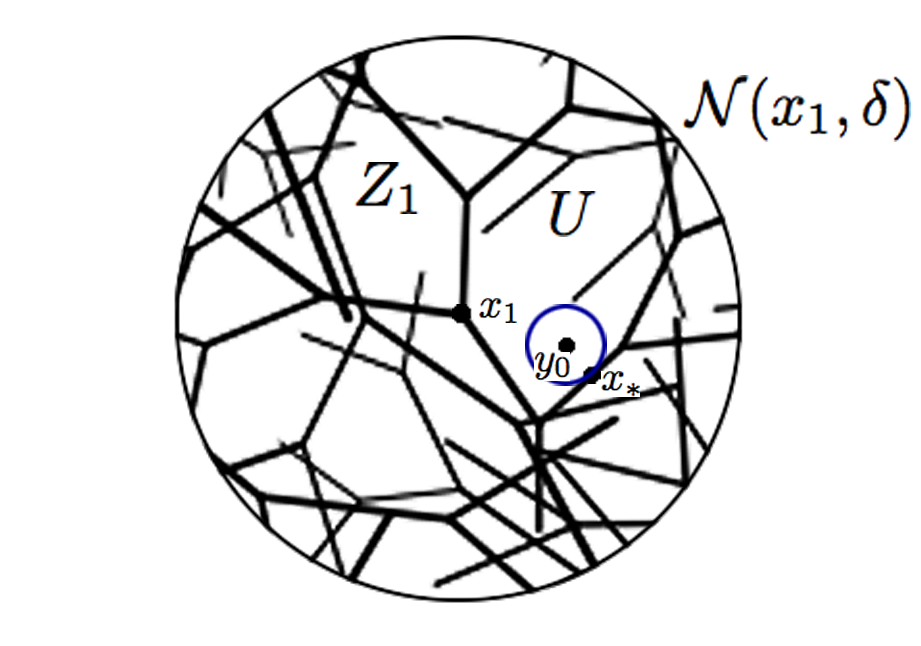}
 \caption{A tame point in the bad set \label{fig:tamepoint}}
\label{Figure1}
\end{center}
\end{figure}

We claim that  $x_*$ is a tame point.
As $\delta >0$ was chosen to be arbitrarily small, and the regular points
are dense in the bad set, the proof of   Proposition~\ref{prop-tame}    then
follows from  this claim.

By the choice of  $B_T(y_0, \delta_0) \subset U$,    there is a   path
$\gamma \colon [0,1] \to  \cN(x_1, \delta)$  such that  $\gamma(0) = y_0$,
$\gamma(1) = x_*$ and $\gamma[0,1) \subset {B_T(y_0, \delta_0)}$.
 The complement of $X_1$ is the good set,
hence  the image  $\gamma[0,1) \subset G$.
The set of leaves with holonomy $G_h$ in $G$ is a union of  submanifolds
with codimension at least $2$ by Proposition~\ref{prop-locstruct}. Thus,  by
a small $C^1$-perturbation of the path $\gamma$ in $U$,
we can assume that its image is disjoint from the set $G_h$. That is,
$\gamma(t) \in   G_e$  for all $0 \leq t < 1$, and $\gamma(1) \in L_*$. Let
$L_t$ denote the leaf containing $\gamma(t)$.

We claim that the volumes of the leaves $L_t$ tend uniformly to infinity.
Assume not, so there exists a constant $V_{max} > 0$ and   a sequence
$0  < t_1 < \cdots <t_n \cdots  \to 1$ such that $x_n = \gamma(t_n) \to x_*$
and the volumes of the leaves
 $L_n = L_{x_n}$  are   bounded above by $V_{max}$. We show this yields a
contradiction to our assumptions.
 What we show in the following is that    if there exists a leaf $L_y$ for
$y \in U \cap G_e$ sufficiently close to $L_*$ with prescribed bounded
volume, then using Proposition~\ref{prop-cptcov}   and
Corollary~\ref{cor-cptcov}, we show this implies that
all leaves intersecting $U$ have bounded volume, which    yields the
contradiction.

\begin{prop}\label{prop-key}
For   $V_{max} > 0$,  there is an $\estar >  0$ so that if there exists $y
\in U \cap G_e$ such that  $d(y, x_*) < \estar $ and $vol(L_y) \leq
V_{max}$,   then for all $y' \in U$, the leaf $L_{y'}$ containing $y'$ has
the volume bound  $vol(L_{y'}) \leq 2  \, V_{max}$.
\end{prop}
 \proof
By Proposition~\ref{prop-volume}, there is a   function $R \colon [0,\infty)
\to [0,\infty)$
  such that if  $L \subset M$ satisfies  $vol(L) \leq V_{max}$ then $\diam
(L) \leq D_* \equiv R(V_{max})$.

   Recall  that $\delta $ was chosen   so  that $2 \delta \leq \etwo$  where
$\e_2 =  \Delta( D_{L_1}, \eone)$ was defined after \eqref{eq-frakp},
   and so that \eqref{eq-choosedelta} holds, hence   $Z_{\F} \subset
\cN(L_1, \eone)$.
   Thus by the choice $x_* \in \overline{B_T(y_0, \delta_0)} \cap  Z_1$, we
have
that $L_* = L_{x_*} \subset  \cN(L_1 , \eone)$.
Let $\pi^* = \pi^{x_*} \colon L_* \to L_1$ denote the normal projection,
whose restriction to $L_*$  is a covering map, which is a diffeomorphism as
$x_* \in Z_1$.

We next choose $y \in U \cap  G_e$ which is sufficiently close to $L_*$ so
that   $L_y \subset  \cN(L_1 , \eone)$ and the holonomy maps of $L_*$ based
at $x_*$ are defined on $L_y$. This will imply that $L_y$ is a finite
covering of $L_*$.

  For each $1 \leq i \leq k$,   let $\tau_i^* \colon [0,1] \to L_*$ be the
lift of  $\tau_i$ with basepoint $x_*$.  By the definition of $D_{L_1}$ in
\eqref{eq-frakp} and the estimate \eqref{eq-unimetric},
  each lifted  path
has bounded length $\| \tau_i^* \| \leq D_{L_1}$ and their homotopy classes
 $\{ [\tau^*_1] , \ldots , [\tau^*_k] \}$ form  a generating set for
$\pi_1(L_* , x_*)$.
 Denote the holonomy along $\tau_i^*$  by $\hh_i^*$.

As $L_*  \subset  \cN(L_1 , \eone)$, there exists
  $0 < \ethree \leq \etwo$ be such that $\cN(L_* , \ethree) \subset \cN(L_1
, \eone)$.

Set $\estar =   \Delta(D_{L_1}, \ethree)$.

By assumption,    there exists $y \in U \cap  G_e \cap  B_T(x_* , \estar) $
with $vol(L_y) \leq V_{max}$, and by
  the choice of $\estar$   we have  $L_y \subset \cN(L_* , \ethree)$.
  Then    the holonomy $\hh_i^*$ along $\tau_i^*$   is   represented by a
map
\begin{equation}\label{eq-extend1}
\hh_i^* \colon \cN(x_*, \estar) \to   \cN(x_*, \ethree) \subset \cN(x_1,
\etwo) \ .
\end{equation}
 Moreover, the bound   $\| \tau_i^* \| \leq D_{L_1}$ implies that   the map
$\hh_i^*$  extends to a  map
 \begin{equation}\label{eq-extend2}
 \hh_i^* \colon \cN(x_*, \etwo) \to    \cN(x_*, \eone) \ .
\end{equation}

 As $L_y$ is a compact leaf, its intersection with the transversal $
\cN(x_*, \ethree)$ is a  finite set, denoted by
\begin{equation}\label{eq-domain}
\fF_y =   L_y \cap  \cN(x_*, \ethree)  \ .
\end{equation}
Then for   each $1 \leq i \leq k$, by \eqref{eq-extend2} the  holonomy   map
$h_{i}^*$ satisfies
 $\ds h_{i}^*  (\fF_y) \subset L_y \cap  \cN(x_*, \ethree) = \fF_y$.
Thus, the finite set of points $\fF_y$ is permuted by the action of a set of
generators for $\pi_1(L_* , x_*)$.
Thus,   compositions of the generators are
defined on  the set
$\fF_y$.
   That is, for any   $w \in \pi_1(L_* , x_*)$ the holonomy $h_w^*$ along
$w$  contains the finite set $\fF_y$ in its domain.
   Let  $\cH_* \subset \pi_1(L_* , x_*)$ denote the normal subgroup of
finite index consisting of all  words whose holonomy    fixes every point in
$\fF_y$.

 Let $z \in \fF_y$.
For each  $w \in \cH_*$,   the holonomy $\hh_w^*$ map is defined at  $z$,
and so must be defined on some open neighborhood $z \in V^w_z \subset U$ of
$z$, where the diameter of the set  $V^w_z$ depends on $z$ and $w$.
As   $y \in U \cap G_e$   the leaf $L_y \subset G_e$ is  without holonomy,
so the restriction of $\hh_w^*$ to the open set $V^w_z$ must fix an open
neighborhood in $\cN(x_*, \eone)$  of $z \in U_z^w \subset V_z^w$. Thus, the
fix-point set of $\hh_w^*$ contains an open neighborhood  of $\fF_y$ in
$\cN(x_*, \eone)$.
Since $y \in L_y \cap  \cN(x_*, \ethree) = \fF_y$, we have
 in particular that there is  an open neighborhood   $y \in U_y^w \subset U
\cap B_T(x_* , \estar)$ contained in the fixed-point set for $\hh_w^*$.

We next use these conclusions for the holonomy of the leaf $L_*$ to deduce
properties of the holonomy for the leaf $L_1$.
Recall that $\cN(L_* , \ethree) \subset \cN(L_1 , \eone)$, and each path
$\tau_i^*$ in $L_*$ is the lift of the path $\tau_i$ in $L_1$ via the
covering map
 $\pi^* \equiv \Pi|L_* \colon L_* \to L_1$. Thus,  the holonomy map
$\hh_i^*$ on $\cN(x_*, \ethree)$ is the restriction of the map $\hh_i$ to
$\cN(x_*, \ethree)$.  Consequently,
the  restriction of  $\hh_i$
to the open set $B_T(x_* , \ethree) \subset \cN(x_1 , \etwo)$ equals the
restriction in \eqref{eq-extend1} of $\hh_i^*$ to $\cN(x_*, \estar)$.
In particular, $\hh_w$ is defined on and fixes  the open set  $U_y^w \subset
B_T(x_* , \estar)$.

 Let $\{w_1 , \ldots, w_N\}$ be a set of generators for $\cH_*$.
Let $m_{\ell}$ denote the word length of $w_{\ell}$ with respect to the
generating set  $\{[\tau_1^*] , \ldots , [\tau_k^*]\}$,
 and set $m_* = \max\{m_1, \ldots ,m_N\}$.

      Fix a choice of $w = w_{\ell} \in \cH_*$.
  Then    the closed path   representing $w$ in  $L_*$ can be lifted to a
path $\tau^y_w$ in   the leaf $L_y$, and as   $L_y \subset   \cN(L_1 ,
\eone)$,
     its  length    is bounded above by       $\|  \tau^y_w\| \leq m_* \cdot
2D_{L_1}$.
   We   show     that  $\hh_w$ is defined on $U$, and $U \subset {\rm
Fix}(\hh_w)$.  This implies that there is a uniform bound on the diameter of
the leaves $L_{y'}$ for $y' \in U$, from which it follows that there is an
upper bound on the function $vol(L_{y'})$ for $y' \in U$, which contradicts
the choice of $U$.

 We first show that  $L_{y'}$   is a finite covering of $L_1$ with the same
index as the covering $L_y \to L_1$.

   Choose $0 < \delta_* \leq  \Delta(2 m_* D_{L_1}, \estar)\leq \estar$ such
that
  $B_T(y, \delta_*) \subset U_y^w$.

 The   open  set $U \subset \cN(x_1,\delta) \smallsetminus Z_1$ is
connected, hence is
path connected. Thus,
given any point $y' \in U$ there is a continuous path $\sigma \colon [0,1]
\to U \cap G_e$ such that $\sigma(0) = y$ and $\sigma(1) = y'$.
  Then   choose a sequence of points $0 = t_0 < t_1 < \cdots < t_{m} = 1$
such that  for $y_i = \sigma(t_i)$, we have   $\sigma([t_i , t_{i+1}])
\subset B_T(y_i , \delta_*)$. See Figure~\ref{fig:goodpoint} below.

\begin{figure}[htbp]
\begin{center}
\includegraphics[width=0.5\textwidth]{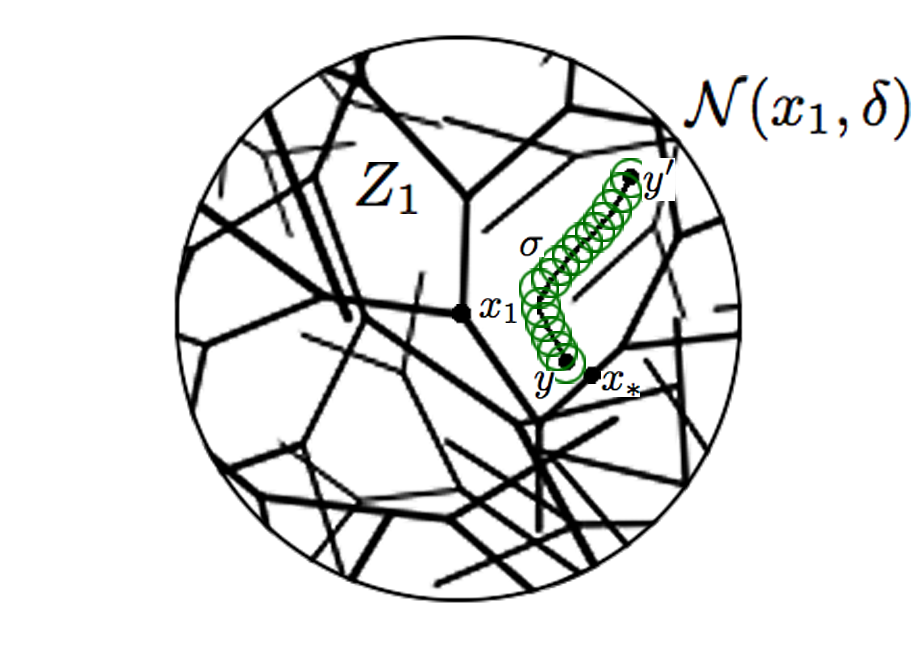}
\caption{A path chain in the good set  \label{fig:goodpoint}}
\label{Figure2}
\end{center}
\end{figure}

We show  that $\sigma([0,1]) \subset {\rm Fix}(\hh_w)$ using induction on
the index $i$. For $i = 0$, $y_0 = y$ and  by assumption, the disk $B_T(y,
\delta_*) \subset U_y^w \subset {\rm Fix}(\hh_w)$ so $\sigma([0,t_1])
\subset {\rm Fix}(\hh_w)$.

Now assume $\sigma([0,t_n])  \subset   {\rm Fix}(\hh_w)$,   hence     $y_n =
\sigma(t_n) \in  {\rm Fix}(\hh_w)$. The closed path $\tau_w^*$ in $L_*$
representing $w$ is the lift of a closed path  $\tau_w$ in $L_1$, which
lifts  to a closed path $\tau_w^{y_n}$ in $L_{y_n}$. As $\tau_w^{y_n}
\subset \cN(L_1 , \eone)$ we have that $\|\tau_w^{y_n}\| \leq 2 m_*
D_{L_1}$.
Then the holonomy map $\hh_w$ for $w$  fixes $y_n$ so near $y_n$ it is
defined by a map
  $$\hh_w^{y_n} \colon \cN(y_n , \delta_*) \to \cN(y_n, \estar) \ . $$

  As the points of $U \cap G_e$ determine leaves without holonomy,
 the set of fixed-points for $\hh_w^{y_n}$  is an open subset of $\cN(y_n ,
\delta_*) \cap U \cap G_e$. The set of fixed-points is also always a
(relatively) closed subset, hence ${\rm Fix}(\hh_w^{y_n})$ contains the
connected component of $\cN(y_n , \delta_*) \cap U \cap G_e$ which contains
the point $y_n$. By assumption we have that
  $\sigma([t_n,t_{n+1}])  \subset  \cN(y_n , \delta_*) \cap U \cap G_e$,
hence
  \begin{equation}
\sigma([t_n,t_{n+1}])  \subset  {\rm Fix}(\hh_w^{y_n}) \subset {\rm
Fix}(\hh_w) \ .
  \end{equation}
 Thus, by induction we conclude that $y' \in  {\rm Fix}(\hh_w)$.

 The choice of $y' \in U$ was arbitrary,   and thus $U \subset  {\rm
Fix}(\hh_w)$.
 We conclude that   $L_{y'}$ is a finite covering of  $L_1$ and isotopic to
$L_y$,  hence $vol(L_{y'}) \leq 2 \, V_{max}$. This completes the proof of
Proposition~\ref{prop-key}.
 \endproof

 The proof of Proposition~\ref{prop-tame} then follows, as the points in the
curve $\sigma(t)$ have limit $x_*$ as $t \to 1$.
 \endproof

\section{Proof of Main Theorem}\label{sec-main}

 In this section, we show that for a compact foliation of a compact
manifold, a categorical open set cannot contain a tame point.
A categorical set must be connected, so we may assume  that $M$ is
connected.
For a connected manifold $M$, there is a finite covering (of degree $d\le
4$) $\wtM \to M$ for which the lifted foliation $\wtF$ is again compact, and
has oriented tangent and normal bundles. We then apply the following two
elementary results to reduce to the oriented case.

\begin{lemma}\label{lem-oriented}
Let $\wtpi \colon \wtM \to M$ be a finite covering with   foliation $\wtF$
whose leaves are finite coverings of the leaves of $\F$. Let $U \subset M$
be a transversely categorical saturated  open  set,
and $H \colon U \times [0,1] \to M$ a foliated homotopy to  the leaf $L_1
\subset M$  of $\F$. Let $\wtU \subset \wtM$ be an open subset such that the
restriction $\wtpi | \wtU \to U$ is a homeomorphism.
Then there exists a foliated homotopy $\wtH \colon \wtU \times [0,1] \to
\wtM$ such that $\wtpi \circ \wtH_t = H_t \circ \wtpi$ for all $0 \leq t
\leq 1$, where $\wtH_t(\wtU) \subset \wtL_1$ for  a finite covering $\wtL_1$
of $L_1$.
\end{lemma}
\proof The covering map $\wtpi$ has the \emph{unique local lifting of paths}
property, so in particular has the \emph{homotopy lifting property}, which
yields the existence of the lifted homotopy $\wtH$.
\endproof

 \begin{lemma}\label{lem-liftedtame}
Let $\wtpi \colon \wtM \to M$ be a finite covering of degree $1 < d <
\infty$, with   foliation $\wtF$ whose leaves are finite coverings of the
leaves of $\F$. Let $\wtL \subset \wtM$ be a leaf of $\wtF$,  and let $\wtx
\in \wtL$ with  $x = \wtpi(\wtx)$. Then $\wtx$ is a tame point in the bad
set for $\wtF$, if and only if  $x$ is a tame point in the bad set for $\F$.
\end{lemma}
\proof
Suppose that $x$ is a tame point for $\F$, then for $\e > 0$  there exists a
continuous path $\gamma \colon [0,1] \to  M$ with $\gamma(1) = x$, as in
Definition~\ref{def-pp}. The map $\wtpi$ has the unique local lifting of
paths property, so there exists a unique path $\wtg \colon [0,1] \to \wtM$
with $\wtg(1) = \wtx$.
Moreover, for each $0 \leq t < 1$ the leaf $\wtL_t$ containing $\wtg(t)$  is
a finite covering of the leaf $L_t \subset M$ containing $\gamma(t)$ where
the degree $\wtpi \colon \wtL_t \to L_t$ has degree at most $d$.  Thus, the
volume $vol(L_t)$ tends to infinity as $t \to 1$, and thus the same holds
for the volume function $\wtvol(\wtL_t)$ in $\wtM$. Thus, $\wtx$ is a tame
point in the bad set for $\wtF$. Conversely, if $\wtx$ is tame point for
$\wtF$ then the proof that $x$ is a tame point for $\F$ follows similarly.
\endproof

 Here is the main result of this section.

 \begin{prop} \label{prop-npp}
Let $\F$ be a compact  $C^1$-foliation of a   compact    manifold $M$. If
$V \subset M$ is a saturated open set which contains a tame point, then $V$
is not transversely categorical.
\end{prop}
\proof
 As a consequence of Lemma~\ref{lem-liftedtame}, we can assume in the
following that both the tangent  bundle $T\F$ and the normal bundle $Q$ to
$\F$ are oriented.

 Let $x_1 \in X_1$ be a tame point, $V\subset  M$ an open set with $x_1 \in
V$, and suppose there exists   a leafwise homotopy $H \colon V \times [0,1]
\to M$ with $H_0 = Id$, and $H_1(V) \subset L_*$ for some leaf $L_*$.
 We    show that this yields  a contradiction.

Recall that for $x \in M$, we let $v(x)$ denote the volume of the leaf $L_x$
containing $x$.

As $x_1$ is a tame point, there is a  smooth  path $\gamma \colon [0,1] \to
V$ such that
$\gamma(1) = x_1$, $\gamma(t) \in G_e$ for $0 \leq t < 1$, and the   volume
$v(\gamma(t))$ of the leaf $L_t$ containing the point $\gamma(t)$
 satisfies $\ds \lim_{t \to 1} v(\gamma(t)) = \infty$.

Define a map  $\phi \colon [0,1] \times [0,1] \to M$ by setting
$\phi_s(t) = \phi(s,t) = H_s(\gamma(t))$.  The key to obtaining a
contradiction is to analyze the behavior of the leaf volume function
$v(\phi(s,t))$.

Set $x_t = \gamma(t)$. Then $x_t \in G_e$ for $0 < t \leq 1$, while $x_1 \in
X_1$ is the given tame point.

As remarked after Definition~\ref{def-pp},    the restricted path
$\gamma \colon [0,1) \to    G_e$
 lies in the set of leaves without holonomy, hence for the leaf $L_0$
containing $x_0 = \gamma(0)$, there is a   foliated isotopy
 $\Gamma  \colon L_{0} \times [0,1)  \to    G_e$ such that $\Gamma_t(x_0) =
x_t$.
In particular,  each map $\Gamma_t \colon L_0 \to L_t$   has homological
degree $1$.

Also note that for $t = 0$, and each $0\leq s \leq 1$, the map   $H_s \colon
L_ 0 =  L_{0,0} \to L_{s,0}$ is surjective  by
Theorem~\ref{thm-compactleaf}.
Let   $d_{s,0}$ denote its homological degree.
  The path of    leaves $s \mapsto L_{s,0}$ starting at $L_0$ has an upper
bound $D_{L_0}$ on their volumes by  Proposition~\ref{prop-bddvol3}, and
moreover,  there is an upper bound
 $ d_{0} = \sup \{ d_{s,0}  \mid 0 \leq s \leq 1\}$.

For $L_1$ the leaf containing the tame point $x_1 = \gamma(1) \in X_1$,
and each $0\leq s \leq 1$,
the map   $H_s \colon L_ 1 =  L_{0,1} \to L_{s,1}$ is also surjective  by
Theorem~\ref{thm-compactleaf}.
Let   $d_{s,1}$ denote its homological degree.
  The path of    leaves $s \mapsto L_{s,1}$ starting at $L_1$ has an upper
bound $D_{L_1}$ on their volumes by  Proposition~\ref{prop-bddvol3}, and
moreover,   there is an upper bound
 $ d_{1} = \sup \{ d_{s,1}  \mid 0 \leq s \leq 1\}$.
   Set
   \begin{equation}\label{eq-volumebounds}
D_L = \max \{D_{L_0} , D_{L_1} \} \ .
\end{equation}

The set $V$ is saturated, so for each  $0 \leq t < 1$, the leaf $L_t \subset
V$ as $\gamma(t) \in V$. Thus, we can define a continuous 2--parameter
family of maps $\Phi \colon [0,1] \times [0,1) \times L_0 \to M$  by setting
 $\Phi_{s,t}(y) =  H_s(\Gamma_t(y))$ for $y \in L_0$.
It is important to recall  the usual caution with the study of compact
foliations:
 the path of leaves $t \mapsto L_{t}$ with unbounded volumes  cannot limit
on  a compact leaf in the bad set.
 Thus, the paths $s \mapsto L_{s, t}$  must become more chaotic   as $t \to
1$, and correspondingly,  the family of maps $\Phi_{s,t}$ is not defined for
$t=1$.
 On the other hand, we are given that the path $\gamma(t)$ limits on $x_1$
and so the trace $\Phi_{s,t}(x_0)$   extends to the continuous map
$\phi(s,t) =  H_s(\gamma(t))$ for $t=1$. We use this   extension  of
$\Phi_{s,t}(y)$ for $y = x_0$ to show that the map extends for all   $y \in
L_0$ which gives a contradiction.

 The idea of the proof of the existence of this extension is to use  the
techniques for studying a homotopy of compact leaves introduced in
Section~\ref{sec-homotopy}, to control the degrees of the maps on the
fundamental classes of the leaves, induced by the maps $\Phi_{s,t}$. This
will in turn yield bounds on the volumes of these leaves, which yields
bounds on their diameters. We can thus use Proposition~\ref{prop-cptcov} to
conclude that for $t_*< 1$ sufficiently close to $t=1$, for each $0 \leq s
\leq 1$, the image   $\Phi_{s,t_*}(L_0)$ is contained in a uniform normal
neighborhood of  $H_s(L_1)$, from which the claim follows. The details
required to fill out this sketch of proof are tedious, but otherwise
straightforward.

First, observe that    $\Phi_{1,t} \colon L_0 \to L_*$,  for $0 \leq t < 1$,
is a family of homotopic maps, hence its homological degree is constant.
Thus, for all  $0 \leq t < 1$ we have:

 \begin{eqnarray*}
\deg(H_1 \colon L_0 \to L_*) & = & \deg(\Phi_{1,0} \colon L_0 \to L_{1,0})\\
 & = &  \deg(\Phi_{1,t} \colon L_0 \to L_{1,t}) \\
 & = &   \deg (\Gamma_t \colon L_0 \to L_t) \cdot \deg(H_1 \colon L_t \to
L_{1,t} = L_*)   \\
 & =&  \deg(H_1 \colon L_t \to L_{1,t})
 \end{eqnarray*}
It follows that
\begin{equation} \label{eq-homolog}
 \deg(H_1 \colon L_t \to L_{1,t}) \leq d_0 ~, ~ \forall ~  0 \leq t < 1 \ .
 \end{equation}
Let $\mathfrak{D} = R(2   d_0  d_1  D_L)$ be the maximum diameter of a leaf
with volume at most
$2 \, d_0 \, d_1 \, D_L$, where we recall that $D_L$ is defined in
\eqref{eq-volumebounds} and $d_0$ and $d_1$ are defined in the text
immediately preceding \eqref{eq-volumebounds}.

 For each $0 \leq s \leq 1$,  recall that  $L_{s,1}$ is the leaf containing
$H_s(x_1)$, and
 let $0 < \e_s' = \e_{L_{s,1}} \leq \ez$ be such that the normal projection
$\Pi_{L_{s,1}}  \colon \cN(L_{s,1}, \e_s') \to L_{s,1}$ is well-defined.
Set $\delta_s'= \Delta(\mathfrak{D} , \e_s')$.

 Let $L$ be a compact  leaf such that  $vol(L) \leq  2 \, d_0 \, d_1 \, D_L$
and  $L \cap \cN(L_{s,1}, \delta_s')  \not= \emptyset$, then by the choice
of $\mathfrak{D}$ and $\delta_s'$,
Proposition~\ref{prop-cptcov}  implies that $L \subset \cN(L_{s,1}, \e_s')$.
  Thus,  the restriction $\Pi_{L_{s,1}} \colon L \to L_{s,1}$ is
well-defined and a covering map, and moreover  by
Corollary~\ref{cor-cptcov}
    we have the estimate
\begin{equation} \label{eq-upper}
vol (L) \leq 2 \deg(\Pi_{L_{s,1}} \colon L \to L_{s,1}) \cdot vol(L_{s,1})
\leq 2 \deg(\Pi_{L_{s,1}} \colon L \to L_{s,1}) \cdot D_L \ .
\end{equation}

 The next step is to choose a   finite covering of the trace of the path
$x_{s,1} = H_s(x_1)$ with respect to the constants $\delta_s'$.
For each $s$, $\cN(L_{s,1}, \delta_s')$ is an open neighborhood of
$L_{s,1}$, so for $\phi(s,t) =  H_s(\gamma(t))$  there
exists $\lambda_s > 0$ such that
\begin{equation}\label{eq-finalcov}
\phi([s-\lambda_s , s + \lambda_s] \times [1-\lambda_s , 1]) \subset
\cN(L_{s,1}, \delta_s') \ .
\end{equation}

Choose a sequence $0 = s_0 <  s_{1} < \cdots < s_{N-1} < s_N  = 1$ of points
such that for $\lambda_n = \lambda_{s_n}$ the collection of open intervals
$\{\cI_n = (s_n - \lambda_n , s_n + \lambda_n) \mid n = 0, 1, \ldots, N\}$
is an open covering of $[0,1]$.

Set $\delta_n'' = \delta_{s_n}'$ and $\e_n'' = \e_{s_n}'$ for $0 \leq n \leq
N$, and
  $\lambda_* = \min\{ \lambda_n \mid n = 0, 1, \ldots, N\} > 0$.

Here is the key result:
\begin{lemma}\label{lem-key}
For $0 \leq s \leq 1$ and  $1 - \lambda_*  \leq t < 1$ we have that
\begin{equation}\label{eq-finalest}
  vol (L_{s,t}) \leq 2 \, d_0 \, d_1   \, D_L \ .
\end{equation}
\end{lemma}
    \proof

For each $1 \leq n \leq N$,   set $\xi_0 =0$ and
  $\xi_{N+1} = 1$, and for $1 \leq n \leq N$ choose    points
$$ \xi_n \in (s_{n-1}  , s_{n-1} + \lambda_{n-1}) \cap (s_{n} - \lambda_{n}
, s_{n})  \ . $$
Then the closed intervals $\{[\xi_0,\xi_1] , [\xi_1,\xi_2], \ldots ,
[\xi_{N-1}, \xi_{N}] , [\xi_{N}, \xi_{N+1}] \}$ form a closed  cover
$[0,1]$.

 Let   $\mu$ satisfy $1 - \lambda_* \leq \mu < 1$, and let $L_{\mu} =
\Gamma_{\mu}(L_0)$ be the leaf through $\gamma(\mu)$.
 The technical idea of the proof of \eqref{eq-finalest}  is to compare the
homological degrees of the maps
 \begin{eqnarray}
H_{\xi_i}|L_{\mu} & \colon & L_{\mu} = L_{0,\mu} \to L_{\xi_i , \mu}
\label{eq-homdegree1} \\
H_{\xi_i}|L_{1} & \colon & L_{1} = L_{0,1} \to L_{\xi_i , 1}
\label{eq-homdegree2}
\end{eqnarray}
using a downward induction argument on $n$, starting with $n = N$, and
showing there is a uniform bound on the ratios of their degrees for all $1 -
\lambda_* \leq \mu < 1$.

For $n =N$,  by \eqref{eq-finalcov}  we have  that
$$\phi([1-\xi_{N}, 1] \times [\mu , 1])  \subset \phi([1-\lambda_{N},1]
\times [1-\lambda_{N} , 1]) \subset \cN(L_{1,1}, \delta_N'')$$
and thus for  each $1-\xi_{N} \leq s \leq 1$
  the point $\phi(s , \mu) \in \cN(L_{1,1}, \delta_N'')$.

   Note that  $L_{1,\mu} = L_{1,1} = L_*$, thus   for $s < 1$ sufficiently
close to $1$ we have $H_s(L_{\mu}) \subset \cN(L_{1,1}, \e_{N}'')$ as the
homotopy $H_s$ is uniformly continuous when restricted to the compact leaf
$L_{\mu}$.

   Let    $r_N$ be the infimum of  $r$  such that     $r \leq s \leq 1$
implies
  $L_{s,\mu} \subset \cN(L_{1,1}, \e_{N}'')$.   The above remark implies
$r_N < 1$.  We  claim  that    $r_N < \xi_{N}$.

Assume, to the contrary, that $r_N \geq \xi_{N}$.
Let $r_N < r < 1$. Then for $ r  \leq s \leq 1$, $L_{s,\mu} \subset
\cN(L_{1,1}, \e_{N}'')$
  and so the normal projection
$\Pi_{L_{1,1}} \colon L_{s,\mu} \to L_{1,1}$ is well-defined and a covering
map.
The restriction
$$H \colon L_{\mu} \times [r, 1] \to  \cN(L_{1,1}, \e_{N}'')$$
yields a homotopy  between
$H_{r} \colon L_{\mu} \to L_{r,\mu}$
and
$H_1 \colon L_{\mu} \to L_{1,\mu} = L_{1,1}$.
Thus,
$$  \deg (\Pi_{L_{1,1}} \circ H_{r} \colon L_{\mu} \to L_{r,\mu} \to
L_{1,1} ) = \deg(\Pi_{L_{1,1}} \circ H_{1} \colon L_{\mu} \to L_{1,\mu} \to
L_{1,1} ) = \deg( H_{1} \colon L_{\mu}  \to L_{1,1} )$$
as
$\Pi_{L_{1,1}} \colon L_{1,\mu} \to L_{1,1}$ is the identity.
The upper bound   \eqref{eq-homolog} implies
$\deg(H_1 \colon L_{\mu} \to  L_{1,1}) \leq d_0$,  hence   the covering
degree of
$\Pi_{L_{1,1}} \colon L_{r,\mu} \to L_{1,1}$  is   bounded above by $d_0$,
as it is an integer which    divides $\deg( H_{1} \colon L_{\mu}  \to
L_{1,1} )$.
  By Corollary~\ref{cor-cptcov}  it follows that
\begin{equation}\label{eq-induct01}
vol (L_{r,\mu}) \leq 2 \, d_0 \cdot vol ( L_{1,1} ) \leq 2 \, d_0\cdot D_L \
.
\end{equation}
The leaf volume function is lower semi-continuous, hence we also have that
$$vol (L_{ r_N, \mu}) ~ \leq ~ \; \lim_{r \to r_N +} ~ vol (L_{r,\mu}) ~
\leq  ~ 2 \, d_0\cdot D_L \ .$$
Thus, the estimate \eqref{eq-induct01} holds for all $r_N \leq r \leq 1$ and
$1-\lambda_* \leq \mu < 1$.

  As we assumed that  $r_N \geq \xi_{N} \geq \lambda_{N}$ we have that
$\phi(r_N , \mu) \in \cN(L_{1,1}, \delta_{N}'')$,
  hence   Proposition~\ref{prop-cptcov}  implies that
  $L_{r_N,\mu} \subset \cN(L_{1,1}, \e_{N}'')$.
  By the uniform continuity of $H_s$ restricted to $L_{\mu}$ at $s = r_N$,
there is
  $r <r_N$ such that   $r < s \leq r_N$ implies
  $L_{s,\mu} \subset \cN(L_{1,1}, \e_{N}'')$. This contradicts the choice of
$r_N$ as the infimum of such $r$, hence we must have that     $r_N <
\xi_{N}$.

  This proves the first statement of the inductive hypothesis for $n=N$,
which is that   the estimate \eqref{eq-induct01}  holds for all $\xi_N  \leq
r \leq 1$ and $1-\lambda_* \leq \mu < 1$.

We next consider    the ratios of covering degrees for a pair of leaves in
adjacent  normal neighborhoods.
 For  $\xi_{N} \leq s \leq 1$,
we have    $\phi(s, 1) \in \cN(L_{1,1}, \delta_{N}'')$ and $vol(L_{s,1})
\leq D_L$
hence  $L_{s,1} \subset \cN(L_{1,1}, \e_{N-1}'')$, and so the normal
projection restricts to a covering map $\Pi_{L_{1,1}} \colon L_{s,1} \to
L_{1,1}$.
Moreover, this implies that  both  $L_{\xi_{N}, \mu}$ and $L_{\xi_{N},1}$
are   coverings of $L_{1,1}$,   and their  homological  degrees are denoted
by
\begin{eqnarray}
\alpha_{N}^{\mu} & = &  \deg (\Pi_{L_{1,1}} \colon L_{\xi_{N},\mu} \to
L_{1,1}) \label{eq-covdegreesN1} \\
a_{N} & = & \deg (\Pi_{L_{1,1}} \colon L_{\xi_{N},1} \to L_{1,1})
\label{eq-covdegreesN2}
\end{eqnarray}

Note that as $s_{N-1} < \xi_N$,
the leaves   $L_{\xi_{N}, \mu}$ and $L_{\xi_{N},1}$ are   also coverings of
$L_{s_{N-1},1}$.
We compare their homological degrees.
 By the uniform continuity of $H_s$ restricted to the curve $\gamma(t)$,
for $0 \leq s \leq 1$, the path $t \mapsto  \phi(s, t)$  has limit $x_{s,1}
= H_{s}(x_1)$.
By Proposition~\ref{prop-cptcov},  the volume bound \eqref{eq-induct01}
for $1-\xi_{N} \leq s \leq 1$  and $1 - \lambda_* \leq t < 1$ implies  that
\begin{equation}\label{eq-path1a}
 H_{s}(L_t) = L_{s, t} \subset  \cN(L_{s,1}, \e_{s}') \ .
 \end{equation}
Thus,  there is a well-defined limit
$$\deg \left( \Phi_{s,1} \colon L_0 \to   L_{s,1} \right) \equiv  \lim_{t\to
1} \left\{ \deg \left( \Pi_{L_{s,1}} \circ \Phi_{s,t} \colon L_0 \to
\cN(L_{s,1}, \e_{s}') \to L_{s,1} \right) \right\} \ .$$
The terminology  $\deg \left( \Phi_{s,1} \colon L_0 \to   L_{s,1} \right)$
is   a small abuse of notation, as
     given  $y \in L_0$ there is no assurance that  $t \mapsto
\Phi_{s,t}(y)$ has a limit at $t=1$; it is only  given  that the image  is
trapped in the open neighborhood $ \cN(L_{s,1}, \e_{s}')$, and the images
are homotopic for $t$ sufficiently close to $1$.

Then for $1-\lambda_s \leq t < 1$, define
\begin{equation}\label{eq-path1b}
 \Xi (s, t)  = \frac{\deg \left(\Phi_{s,1} \colon L_0 \to L_{s,1}
\right)}{\deg \left( \Phi_{s, t} \colon L_0 \to L_{s,t} \right)} \ .
  \end{equation}

 We now apply this  discussion in the case   $s = \xi_{N}$ where we have the
volume bound \eqref{eq-induct01}.
  It  again follows from Proposition~\ref{prop-cptcov}    that   for
$1-\lambda_* \leq t < 1$,  and noting that $s_N =1$,
\begin{equation}\label{eq-path1c}
 H_{\xi_{N}}(L_t) = L_{\xi_{N}, t} \subset  \cN(L_{s_N,1}, \e_{N}'') \cap
\cN(L_{{s_{N-1}},1}, \e_{N-1}'') \ .
 \end{equation}
Thus, for $1-\lambda_* \leq \mu \leq  t < 1$ the maps
$$\Pi_{L_{1,1}} \circ \Phi_{\xi_{N},\mu} \sim \Pi_{L_{1,1}} \circ
\Phi_{\xi_{N},t}  \colon L_0 \to \cN(L_{1,1}, \e_{N}'')$$
are homotopic, hence
\begin{equation}\label{eq-induct1c}
\deg ( \Pi_{L_{1,1}} \circ \Phi_{\xi_{N},\mu}) = \deg ( \Pi_{L_{1,1}} \circ
\Phi_{\xi_{N},t}) \ .
\end{equation}
 For $t$ sufficiently close to $1$ the  map $ \Pi_{L_{1,1}} \circ
\Phi_{\xi_{N},t}$ on the left-hand-side of \eqref{eq-induct1c}   factors
$$\Pi_{L_{1,1}} \circ \iota \circ \Pi_{L_{\xi_{N},1}} \circ \Phi_{\xi_{N},t}
\colon L_0 \to \cN(L_{{\xi_{N}},1}, \e_{\xi_{N-1}}')  \to L_{\xi_{N},1}
\subset \cN(L_{1,1}, \e_{N}'') \to L_{1,1}$$
while the  map $\ds \Pi_{L_{1,1}} \circ \Phi_{\xi_{N},\mu} $ on
right-hand-side of \eqref{eq-induct1c}   factors
$$\Pi_{L_{1,1}} \circ \Phi_{\xi_{N},\mu} \colon L_0 \to L_{{\xi_{N}},\mu}
\to L_{1,1} \ .$$
Identifying the degrees of these maps in our terminology, we obtain from
\eqref{eq-induct1c}  that
$$ \deg( \Phi_{\xi_N,\mu} \colon L_0 \to L_{{\xi_{N}},\mu}) \cdot
\alpha_{N}^{\mu}  = \deg ( \Pi_{L_{1,1}} \circ \Phi_{\xi_{N},\mu}) = \deg
( \Pi_{L_{1,1}} \circ \Phi_{\xi_{N},t})  =  \deg( \Phi_{\xi_N,1} \colon L_0
\to L_{{\xi_{N}},1}) \cdot a_{N}$$
and so
\begin{equation}\label{eq-induct-xi1}
\alpha_{N}^{\mu} = \Xi(\xi_{N}, \mu) \cdot  a_{N} \ .
\end{equation}
Thus, the ratio \eqref{eq-path1b} gives the relation between the homological
degrees of the maps in \eqref{eq-covdegreesN1} and \eqref{eq-covdegreesN2}.
 This completes the proof of the first stage of the induction.

The general    inductive hypotheses involves two statements: Given  $0 \leq
n \leq N$, we first  assume that:
 \begin{equation}
{\rm for ~ all} ~ 0 \leq n \leq N,  ~ {\rm for ~ all} ~ \xi_{n} \leq s \leq
1  ~ {\rm and} ~  1 - \lambda_* \leq t \leq 1  , ~ {\rm then} ~ vol
(L_{s,t})    \leq    2 \, d_0 \, d_1\cdot D_L \ .  \label{eq-inductn1}
\end{equation}
  Given  \eqref{eq-inductn1},  then for   $n \leq \ell \leq N$ and $1 -
\lambda_* \leq \mu \leq 1$ define  the     integers $a_{\ell} , b_{\ell}  ,
\alpha_{\ell}^{\mu} , \beta_{\ell}^{\mu}$.
$$
 \begin{array}{cclclcl}
 L_{\xi_{\ell}, 1} & \subset & \cN(L_{s_{\ell},1}, \e_{\ell}'') & , &
a_{\ell} & = & \deg \left(\Pi_{L_{s_{\ell},1}} \colon  L_{\xi_{\ell}, 1} \to
L_{s_{\ell},1} \right) \\
  L_{\xi_{\ell}, 1} & \subset & \cN(L_{s_{\ell -1},1}, \e_{\ell -1}'')  & ,
&  b_{\ell} & = &  \deg \left( \Pi_{L_{s_{\ell -1},1}} \colon
L_{\xi_{\ell}, 1} \to L_{s_{\ell -1},1}   \right) \\
 L_{\xi_{\ell}, \mu} & \subset & \cN(L_{s_{\ell},1}, \e_{\ell}'') & , &
\alpha_{\ell}^{\mu}  & = &   \deg \left( \Pi_{L_{s_{\ell},1}} \colon
L_{\xi_{\ell}, \mu}  \to  L_{s_{\ell},1}    \right) \\
L_{\xi_{\ell}, \mu} & \subset & \cN(L_{s_{\ell -1},1}, \e_{\ell -1}'') & , &
\beta_{\ell}^{\mu} & = &    \deg \left( \Pi_{L_{s_{\ell -1},1}} \colon
L_{\xi_{\ell}, \mu}  \to L_{s_{\ell -1},1}   \right)
 \end{array}
 $$
For notational convenience, set $b_{N+1} = \beta_{N+1}^{\mu} = 1$ and $a_{0}
= \alpha_{0}^{\mu} = 1$. Second, we assume that:
   \begin{equation}\label{eq-inductn2}
{\rm for ~ all} ~  n \leq \ell \leq N, ~ {\rm and} ~ 1 - \lambda_* \leq \mu
\leq 1, ~ {\rm then} ~
\frac{\alpha_{\ell}^{\mu}}{a_{\ell}}  =   \Xi(\xi_{\ell}, \mu) = \frac{
\beta_{\ell}^{\mu}}{b_{\ell}} \ .
\end{equation}
We   show that  if \eqref{eq-inductn1}  and \eqref{eq-inductn2}  are true
for $n$, then the corresponding statements   are true   for $n-1$.

The choice of $\lambda_s > 0$ so that      \eqref{eq-finalcov}  holds
implies   that
$$\phi([s_{n-1}-\lambda_{n-1}, s_{n-1}+\lambda_{n-1}] \times [1-\lambda_* ,
1])    \subset \cN(L_{s_{n-1},1}, \delta_{n-1}'')$$
and hence    $\phi(s, t) \in \cN(L_{s_{n-1},1}, \delta_{n-1}'')$ for all
$\xi_{n-1} \leq s \leq \xi_{n}$ and $1-\lambda_* \leq t < 1$.

For $s = \xi_n$ the hypothesis  \eqref{eq-inductn1}  implies that  for all
$1 - \lambda_* \leq t < 1$,
\begin{equation}
 vol (L_{\xi_n,t})    \leq    2 \, d_0 \, d_1\cdot D_L ~ \; {\rm and} ~ {\rm
hence}  \; ~  L_{\xi_n,t}    \subset    \cN(L_{s_{n-1},1}, \e_{n-1}'') \ .
 \end{equation}
Thus, the restriction $\ds \Pi_{L_{s_{n-1},1}} \colon  L_{\xi_n,t} \to
L_{s_{n-1},1}$ is a covering map. The key to the proof of the inductive step
is to obtain a uniform estimate for the homological degree of this covering
map.

\begin{lemma} \label{lem-upperbd}
For all ~ $1 - \lambda_* \leq t < 1$,
$\ds \beta_n^t \cdot   \deg \left( H_{\xi_{n}} \colon L_{0,t} \to L_{\xi_{n}
, t} \right) ~  \leq  ~   d_0 \, d_1$.
\end{lemma}
\proof  Consider  the diagram

\begin{picture}(600,100)\label{eqn.diag}

\put(10,80){$L_{0,t}$}

\put(40,82){\vector(1,0){100}}
\put(90,90){$H_{\xi_n}$}

\put(160,80){$L_{\xi_n,t}$}

\put(280,80){$L_{\xi_{n+1},t}$}

\put(335,80){$\cdots$}

\put(370,80){$L_{\xi_{N},t}$}

\put(10,20){$L_{0,1}$}

\put(40,22){\vector(1,0){40}}
\put(48,10){$H_{s_{n-1}}$}

\put(90,20){$L_{s_{n-1},1}$}

\put(150,22){\vector(-1,0){20}}
\put(138,10){$b_n$}

\put(160,20){$L_{\xi_n,1}$}

\put(190,22){\vector(1,0){20}}
\put(196,10){$a_n$}

\put(220,20){$L_{s_n,1}$}

\put(270,22){\vector(-1,0){20}}
\put(257,10){$b_{n+1}$}

\put(280,20){$L_{\xi_{n+1},1}$}

\put(335,20){$\cdots$}

\put(370,20){$L_{\xi_{N},1}$}

\put(410,22){\vector(1,0){20}}
\put(416,10){$a_N$}

\put(440,20){$L_{1,1}$}

\put(152,72){\vector(-1,-1){37}}

\put(167,62){$\vdots$}
\put(167,50){$\vdots$}
\put(168,47){\vector(0,-1){12}}

\put(182,72){\vector(1,-1){37}}

\put(280,72){\vector(-1,-1){37}}

\put(294,62){$\vdots$}
\put(294,50){$\vdots$}
\put(295,47){\vector(0,-1){12}}

\put(400,72){\vector(1,-1){37}}

\put(379,62){$\vdots$}
\put(379,50){$\vdots$}
\put(380,47){\vector(0,-1){12}}

\put(115,55){$\beta_{n}^t$}

\put(170,40){$\Xi (n,t)$}

\put(205,55){$\alpha_{n}^t$}

\put(235,55){$\beta_{n+1}^t$}

\put(297,40){$\Xi (n+1,t)$}

\put(382,40){$\Xi (N,t)$}

\put(425,55){$\alpha_{N}^t$}

\end{picture}

where the integer next to a covering map indicates its homological degree.

The maps $H_{\xi_n} \colon L_{0,1} \to L_{\xi_n , 1}$ and
$H_{s_{n-1}} \colon L_{0,1} \to L_{s_{n-1} , 1}$ are homotopic through maps
into $\cN(L_{s_{n-1} , 1}, \e_{n-1}'')$,  hence
\begin{equation} \label{eq-71}
\deg \left( H_{s_{n-1}} \colon L_{0,1} \to L_{s_{n-1} , 1} \right) = b_n
\cdot \deg \left( H_{\xi_n} \colon L_{0,1} \to L_{\xi_n , 1} \right) \ .
\end{equation}
As $\ds \deg \left( H_{s_{n-1}} \colon L_{0,1} \to L_{s_{n-1} , 1} \right) =
d_{s,1}  \leq d_1$ and the degrees of the maps are positive integers, it
follows that $1 \leq b_n \leq d_1$.

For $n \leq \ell < N$ and  $1 - \lambda_* \leq t < 1$,
the maps $H_{\xi_{\ell}} \colon L_{0,t} \to L_{\xi_{\ell} , t}$ and
$H_{\xi_{\ell + 1}} \colon L_{0,t} \to L_{\xi_{\ell +1} , t}$ are homotopic
through maps into $\cN(L_{s_{\ell} , 1}, \e_{\ell}'')$,  hence
\begin{equation}  \label{eq-72}
\alpha_{\ell}^t \cdot \deg \left( H_{\xi_{\ell}} \colon L_{0,t} \to
L_{\xi_{\ell} , t} \right) = \beta_{\ell +1}^t  \cdot \deg \left(
H_{\xi_{\ell +1}} \colon L_{0,t} \to L_{\xi_{\ell +1} , t} \right) \ .
\end{equation}

Likewise, for $n \leq \ell < N$,
the maps $H_{\xi_{\ell}} \colon L_{0,1} \to L_{\xi_{\ell} , 1}$ and
$H_{\xi_{\ell + 1}} \colon L_{0,1} \to L_{\xi_{\ell +1} , 1}$ are homotopic
through maps into $\cN(L_{s_{\ell} , 1}, \e_{\ell}'')$,  hence
\begin{equation}  \label{eq-73}
a_{\ell} \cdot \deg \left( H_{\xi_{\ell}} \colon L_{0,1} \to L_{\xi_{\ell} ,
1} \right) = b_{\ell +1}  \cdot \deg \left( H_{\xi_{\ell +1}} \colon L_{0,1}
\to L_{\xi_{\ell +1} , 1} \right) \ .
\end{equation}

It follows from equation   \eqref{eq-72}     that
\begin{eqnarray*}
\deg \left( H_{1} \colon L_{0,t} \to L_{1 , t} \right)  & = &
\frac{ \alpha_{N}^t}{\beta_{N +1}^t} \cdot \deg \left( H_{\xi_{N}} \colon
L_{0,t} \to L_{\xi_{N} , t} \right)\\
  & = &
\frac{ \alpha_{N-1}^t \alpha_{N}^t}{ \beta_{N}^t \beta_{N +1}^t} \cdot \deg
\left( H_{\xi_{N-1}} \colon L_{0,t} \to L_{\xi_{N-1} , t} \right)\\
& \vdots  & \\
  & = &
\frac{\alpha_n^t \cdots  \alpha_{N-1}^t \alpha_{N}^t}{\beta_{n+1}^t \cdots
\beta_{N}^t \beta_{N +1}^t} \cdot \deg \left( H_{\xi_{n}} \colon L_{0,t} \to
L_{\xi_{n} , t} \right) \\
  & = &
\frac{\alpha_{n}^t \cdots  \alpha_{N}^t}{\beta_{n}^t \cdots  \beta_{N}^t}
\cdot \beta_n^t \cdot   \deg \left( H_{\xi_{n}} \colon L_{0,t} \to
L_{\xi_{n} , t} \right)
\end{eqnarray*}
so that by the inductive hypothesis \eqref{eq-inductn2}  we have
\begin{eqnarray}
\beta_n^t \cdot   \deg \left( H_{\xi_{n}} \colon L_{0,t} \to L_{\xi_{n} , t}
\right)  & =  &
\frac{\beta_{n}^t \cdots  \beta_{N}^t}{\alpha_{n}^t \cdots  \alpha_{N}^t}
\cdot \deg \left( H_{1} \colon L_{0,t} \to L_{1 , t} \right) \\
   & =  &
 \frac{b_n \cdots b_N}{a_n \cdots a_N}   \cdot \deg \left( H_{1} \colon
L_{0,t} \to L_{1 , t} \right) \label{eq-74}
\end{eqnarray}

Using   \eqref{eq-73}  we obtain

\begin{equation} \label{eq-75}
  \deg \left( H_{1} \colon L_{0,1} \to L_{1 , 1} \right) = \frac{a_n \cdots
a_N}{b_n \cdots b_N}  \cdot \deg \left( H_{s_{n-1}} \colon L_{0,1} \to
L_{s_{n-1} , 1}  \right) \ .
  \end{equation}

so that

\begin{equation} \label{eq-76}
\frac{b_n \cdots b_N}{a_n \cdots a_N}   =  \frac{\deg \left( H_{s_{n-1}}
\colon L_{0,1} \to L_{s_{n-1} , 1}  \right)}{ \deg \left( H_{1} \colon
L_{0,1} \to L_{1 , 1} \right)} \leq d_1 \ .
\end{equation}

and hence combining \eqref{eq-homolog} ,  \eqref{eq-74}  and \eqref{eq-76}
we obtain

\begin{equation} \label{eq-77}
\beta_n^t \cdot   \deg \left( H_{\xi_{n}} \colon L_{0,t} \to L_{\xi_{n} , t}
\right)   \leq  d_1  \cdot  \deg \left( H_{1} \colon L_{0,t} \to L_{1 , t}
\right)  \leq d_0 \, d_1 \ .
\end{equation}
This completes the proof of Lemma~\ref{lem-upperbd}.
\endproof

 Fix $1 - \lambda_* \leq \mu < 1$.
Let $r_{n-1} \leq \xi_n$ be the infimum  of $r$ satisfying   $  r \leq
\xi_n$ such that
  $r \leq s \leq \xi_n$ implies that
  $L_{s,\mu}  \subset     \cN(L_{s_{n-1},1}, \e_{n-1}'')$.
As  $ L_{\xi_n,\mu}    \subset    \cN(L_{s_{n-1},1}, \e_{n-1}'')$, the
continuity of $H_s$ at $s = \xi_n$ implies that
 $r_{n-1} < \xi_n$. We claim that $r_{n-1} < \xi_{n-1}$.

Assume, to the contrary, that $r_{n-1} \geq \xi_{n-1}$.
Let $r_{n-1} < r < \xi_n$, then for $ r  \leq s \leq \xi_n$,
$L_{s,\mu}\subset     \cN(L_{s_{n-1},1}, \e_{n-1}'')$
  and so the normal projection
$\Pi_{L_{s_{n-1},1}} \colon L_{s,\mu} \to L_{s_{n-1},1}$ is well-defined and
a covering map.
The restriction
$$H \colon L_{\mu} \times [r, \xi_n] \to   \cN(L_{s_{n-1},1}, \e_{n-1}'')$$
yields a homotopy  between
$H_{r} \colon L_{\mu} \to L_{r,\mu}$
and
$H_{\xi_n} \colon L_{\mu} \to L_{\xi_n,\mu}$.
Thus,
$$  \deg (\Pi_{L_{s_{n-1},1}} \circ H_{r} \colon L_{\mu} \to L_{r,\mu} \to
L_{s_{n-1},1} ) = \deg(\Pi_{L_{\xi_{n-1},1}} \circ H_{\xi_n} \colon L_{\mu}
\to L_{\xi_n,\mu} \to L_{\xi_{n-1},1} ) \ . $$

It follows from the estimate \eqref{eq-77}  that
\begin{equation}\label{eq-78}
\deg (\Pi_{L_{s_{n-1},1}} \colon L_{r,\mu} \to L_{s_{n-1},1} ) \leq \deg (
\Pi_{L_{s_{n-1},1}} \circ H_{r} \colon L_{\mu} \to L_{r,\mu} \to
L_{s_{n-1},1}) \leq d_0 \, d_1
\end{equation}
hence
\begin{equation}\label{eq-79}
vol ( L_{r,\mu})  \leq 2 \, d_0 \, d_1 \cdot vol ( L_{s_{n-1},1}) \leq  2
\,d_0 \, d_1 \cdot D_L \ .
\end{equation}

The leaf volume function is lower semi-continuous, hence we also have that
\begin{equation}\label{eq-80}
vol (L_{ r_{n-1}, \mu}) ~ \leq ~ \; \lim_{r \to r_{n-1} +} ~ vol (L_{r,\mu})
~ \leq  ~ 2 \, d_0 \, d_1 \cdot D_L \ .
\end{equation}
Thus, the estimate \eqref{eq-79}  holds for all $r_{n-1} \leq r \leq 1$ and
$1-\lambda_* \leq \mu < 1$.

  As we assumed that  $r_{n-1} \geq \xi_{n-1} \geq s_{n-1} - \lambda_{n-1}$
we have that  $\phi(r_{n-1} , \mu) \in \cN(L_{s_{n-1},1}, \delta_{n-1}'')$
  hence
  $L_{r_{n-1},\mu} \subset \cN(L_{s_{n-1},1}, \e_{n-1}'')$.
  By the continuity of $H_s$ at $s = r_{n-1}$ there is
  $r <r_{n-1}$ such that   $r < s \leq r_{n-1}$ implies
  $L_{s,\mu} \subset \cN(L_{s_{n-1},1}, \e_{n-1}'')$. This contradicts the
choice of $r_{n-1}$ as the infimum of such $r$, hence we must have that
$r_{n-1} < \xi_{n-1}$.  This proves the first statement of the inductive
hypothesis for $n-1$.

The second inductive statement \eqref{eq-inductn2}  follows exactly as
before.

Thus, we conclude by downward induction that   \eqref{eq-finalest}  holds
for all $1-\lambda_* \leq t < 1$ and all $0 \leq s \leq 1$.
\endproof

The conclusion \eqref{eq-finalest} of Lemma~\ref{lem-key} for $s = 0$
contradicts the assumption  that    $\ds \lim_{t \to \infty} \ vol(L_{0,t})
= \infty$.
This   contradiction  completes the proof of Proposition~\ref{prop-npp}.
\endproof

\vfill
\eject


\end{document}